\documentclass[12pt,a4paper,draft]{article}
\usepackage{latexsym,amsfonts,amssymb,amsmath,longtable,amsthm}
\renewcommand{\le}{\leqslant}
\renewcommand{\ge}{\geqslant}

\begin{document}

\newtheorem{lem}{{\bfseries Lemma}}[section]
\newtheorem{ttt}[lem]{{\bfseries Theorem}}
\newtheorem{cor}[lem]{{\bfseries Corollary}}
\newtheorem{prp}[lem]{{\bfseries Proposition}}
\theoremstyle{remark}
\newtheorem*{rem}{{\bfseries Remark}}
\newtheorem*{quest}{{\bfseries Problem}}
\theoremstyle{definition}
\newtheorem{tab}[lem]{{\bfseries Table}}
\newtheorem{df}[lem]{\bfseries Definition}
\newtheorem{con}[lem]{\slshape Conjecture}

\newcommand{\Sym}{\mathrm{Sym}}
\newcommand{\F}{\mathbb{F}}
\newcommand{\Z}{\mathbb{Z}}
\newcommand{\ov}{\overline}
\newcommand{\Aut}{\mathrm{Aut} }
\renewcommand{\P}{\mathbf{P}}

\begin{flushright}
UDC 512.542
\end{flushright}

\begin{center}
{\bfseries \Large Carter subgroups of finite almost simple groups}

E.P.Vdovin
\end{center}

\begin{center}
{\bfseries Annotation}
\end{center}

In the paper we complete the classification of Carter subgroups in finite almost simple
groups. In particular, we prove that Carter subgroups of every finite almost simple group
are conjugate. Togeather with previous results by author and F. Dalla Volta, A. Lucchini,
and M. C. Tamburini, as a corollary, it follows that Carter subgroups of every finite
group are conjugate.

\section{Introduction}

We recall that a subgroup of a finite group is called a {\em Carter subgroup} if it
is nilpotent and self-normalizing. By a well-known result, any finite solvable group
contains exactly one conjugacy class of Carter subgroups (cf. \cite{Ca1}), and it is
reasonable to conjecture that a finite group contains at most one conjugacy class of
Carter subgroups. The evidence for this conjecture is based on extensive
investigation, by several authors, of classes of finite groups which are close to be
simple. In particular it has been shown that the conjecture holds for the symmetric
and alternating groups (cf. \cite{DT1}) and, denoting by $p^t$ a power of a prime
$p$, for any group $A$ such that $SL_n(p^t)\leq A\leq GL_n(p^t)$ (cf. \cite{DT2} and
\cite{Va}), for the symplectic groups $Sp_{2n}(p^t)$, the full unitary groups
$GU_n(p^{2t})$ and, when $p$ is odd, the full orthogonal groups $GO_n^{\pm}(p^t)$
(cf. \cite{DTZ}). Later in \cite{PreTamVdo} results of \cite{DTZ} were extended to
any group $G$ with $O^{p'}(S)\leq G\leq S$, where $S$ is a full classical matrix
group. Also some of the sporadic simple groups were investigated (cf.  \cite{D'A},
for example). In the nonsolvable cases, when Carter subgroups exist, they always
turn out to be the normalizers of Sylow $2$-subgroups.

In the paper we consider the following

\begin{quest}
Are any two Carter subgroups of a finite group conjugate?
\end{quest}

In \cite{DLT} it is proven that the minimal counterexample $A$ to this
problem should be almost simple. Later in \cite{Vdo} a stronger result was
obtained. A finite group $G$ is said to satisfy condition
($*$) if, for every its non-Abelian composition factor $S$ and
for every its nilpotent subgroup $N$, Carter subgroups of $\langle
\mathrm{Aut}_N(S),S\rangle$ are conjugate (definition of $\Aut_N(S)$ one can
find below). In \cite{Vdo} the following theorem was proven.

\begin{ttt}\label{maininduct}
If a finite group $G$ satisfies {\em ($*$)}, then Carter subgroups of
$G$ are conjugate.
\end{ttt}

 Thus our goal here is to prove that for every known simple group $S$ and every
nilpotent subgroup $N$ of $\Aut(S)$, Carter subgroups of $\langle S,N\rangle$
are conjugate.  Some classes
of almost simple groups which can not be minimal counter example to the problem
are found in \cite{PreTamVdo} and~\cite{TamVdo}. The resulting table of almost
simple groups with conjugate Carter subgroups is given in~\cite{Vdo}.

Our notations is standard.  If $G$ is a finite group, we denote  by $\P G$
the factor group $G/Z(G)$. If $\pi$ is a set of primes then we denote by $\pi'$
its complement in the set of all primes. As usual we denote by $O_\pi(G)$ the
maximal normal $\pi$-subgroup of $G$ and we denote by $O^{\pi'}(G)$ the
subgroup generated by all $\pi$-elements of $G$. If $\pi=\{2\}'$ is the set of
all odd primes, then $O_\pi(G)=O_{2'}(G)$  is denoted by~$O(G)$. If $g\in G$,
 then we denote by $g_\pi$ the $\pi$-part of $g$, i.~e., $g_\pi=g^{\vert
g\vert_{\pi'}}$. For a finite
group $G$ we denote by $\Aut(G)$ the
group of automorphisms of $G$. If $\lambda\in\Aut(G)$, then we
denote by $G_\lambda$ the set of $\lambda$-stable points, i.~e.,
$G_\lambda=\{g\in
G\vert g^\lambda=g\}$. If $Z(G)$ is trivial, then $G$ is
isomorphic to the group of its inner automorphisms and we may suppose
that~$G\leq \Aut(G)$.  A finite group $G$ is said to be {\em almost
simple} if there is a simple group $S$ with $S\leq G\leq \Aut(S)$, i.~e.,
$F^\ast(G)$ is a simple group. We denote by $F(G)$ the Fitting subgroup of $G$
and by $F^\ast(G)$ the generalized Fitting subgroup of~$G$.

If $G$ is a group, $A,B,H$ are subgroups of $G$ and $B$ is normal in $A$
($B\unlhd A$), then $N_H(A/B)=N_H(A)\cap N_H(B)$. If $x\in N_H(A/B)$, then $x$
induces an automorphism $Ba\mapsto B x^{-1}ax$ of $A/B$. Thus, there is a
homomorphism of $N_H(A/B)$ into $\Aut(A/B)$. The image of this homomorphism is
denoted by $\Aut_H(A/B)$ while its kernel is denoted by $C_H(A/B)$. In
particular, if $S$ is a composition factor of $G$, then for any $H\leq G$ the
group $\Aut_H(S)$ is defined.

\section{Preliminary results}

\begin{lem}\label{HomImageOfCarter}
Let $G$ be a finite group, let  $K$ be a Carter subgroup of $G$ and assume
that $N$ is a normal subgroup of $G$. Assume that $KN$ satisfies {\em($*$)}.
Then $KN/N$ is a Carter subgroup of~$G/N$.
\end{lem}

\begin{proof}
Consider $x\in G$ and assume that $xN\leq N_{G/N}(KN/N)$. It follows
that $x\in N_G(KN)$.  We have that
$K^x$ is a Carter subgroup of $KN$. Since $KN$ satisfies ($*$), we have
that its Carter subgroups are
conjugate. Thus there exists $y\in KN$ such that $K^y=K^x$. Since $K$ is a
Carter subgroup of $G$, it follows that $xy^{-1}\in N_G(K)=K$ and~${x\in KN}$.
\end{proof}

\begin{lem}\label{power} {\em \cite[Lemma~5]{Vdo}}
Assume that $G$ is a finite group. Let $K$ be a Carter subgroup of
$G$, with centre $Z(K)$. Assume also that  $e\not= z\in Z(K)$ and $C_G(z)$
satisfies~{\em($*$)}.
\begin{itemize}
\item[{\em (1)}] Every subgroup $Y$ which contains $K$ and satisfies {\em ($*$)}
is self-normalizing in~$G$.
\item[{\em (2)}] No conjugate of $z$ in $G$, except $z$,
lies in $Z(K)$.
\item[{\em (3)}] If $H$ is a Carter subgroup of $G$, non-conjugate
to $K$, then $z$ is not conjugate to any
element in the centre of $H$.
\end{itemize}

In particular the centralizer $C_G(z)$ is self-normalizing
in $G$, and $z$ is not conjugate to any power~${z^k\not= z}$.
\end{lem}

\begin{lem}\label{CritSyl2Carter}
Let $G$ be a finite group and $S$ be a Sylow $2$-subgroup of $G$. Then
$G$ contains a Carter subgroup $K$ with $S\leq K$ if and only
if~$N_G(S)=SC_G(S)$.
\end{lem}

\begin{proof}
Assume that $G$ contains a Carter subgroup $K$ with $S\leq K$. Since  $K$ is
nilpotent, it follows that $S$ is normal in $K$ and $K\leq SC_G(S)\unlhd N_G(S)$. By
Feit-Thompson Theorem (see \cite{ft}) we obtain that $N_G(S)$ is solvable. Thus, by
Lemma \ref{power}(1) we have that $SC_G(S)$ is self-normalizing in $G$,
therefore~${N_G(S)=S C_G(S)}$.

Assume now that $N_G(S)=SC_G(S)$, i.~e., $N_G(S)=S\times
O(C_G(S))$. Since $O(C_G(S))$ is of odd order, it is
solvable. Hence it contains a Carter subgroup $K_1$. Consider a
nilpotent subgroup $K=S\times K_1$ of $G$. Assume that $x\in
N_G(K)$, then $x\in N_G(S)$. But $K$ is a Carter subgroup of $N_G(S)$,
hence $x\in K$ and $K$ is a Carter subgroup of~$G$.
\end{proof}

\begin{lem}\label{Syl2centrcomposit}
Let $G$ be a finite group, let $S$ be a Sylow $2$-subgroup of $G$ and $x\in
N_G(S)$
be of odd order. Assume that there exist normal subgroups $G_1,\ldots,G_k$ of
$G$
such that $G_1\cap\ldots\cap G_k\cap S\leq Z(N_G(S))$. If
$\varphi_i:G\rightarrow G_i$ is the natural homomorphism assume also that
$x^{\varphi_i}$ centralizes~$SG_i/G_i$.

Then $x$ centralizes $S$.
\end{lem}

\begin{proof}
Consider the normal series $S\rhd S_1\rhd \ldots\rhd S_k\rhd S_{k+1}=\{e\}$, where
$S_i=S\cap(G_1\cap\ldots\cap G_i)$. The conditions of the lemma implies that $x$
centralizes every factor $S_{i-1}/S_i$. Since $x$ is of odd order this imply that
$x$ centralizes~$S$.
\end{proof}

\begin{lem}\label{CarterInDirectProduct} {\em \cite[Lemma~3]{Vdo}}
Let $G$ be a finite group. Let $H$ be a Carter
subgroup of $G$. Assume that there exists a normal subgroup
$B=T_1\times\ldots\times T_k$ of $G$ such that $T_1\simeq\ldots\simeq T_k\simeq
T$, $Z(T_i)=\{1\}$ for all $i$, and $G=H(T_1\times\ldots\times T_k)$. Then
$\mathrm{Aut}_H(T_i)$ is a Carter subgroup of~$\langle \mathrm{Aut}_H(T_i),
T_i\rangle$.
\end{lem}

\begin{lem}\label{InhBy2-ext}
Let $G$ be a finite group, let $H$ be a normal subgroup of $G$  such that
$\vert G:H\vert=2^t$. Let $S,T$ be Sylow $2$-subgroups of $G,H$ respectively
and $N_H(T)=TC_H(T)$. Then~$N_G(S)=SC_G(S)$.

In particular, both $G,H$ contain Carter subgroups $K,L$ respectively with
$S\leq K$ and~${T\leq L}$.
\end{lem}

\begin{proof}
Consider $N_G(S)$. Since $H$ is normal in $G$ we have
that $$N_G(S)\leq N_G(T)=\langle S,N_H(T)\rangle=\langle S,T\times
O(N_H(T))\rangle.$$ Since $N_H(T)$ is normal in $N_G(T)$, we have that
$O(N_H(S))=O(C_H(T))$ is normal in $N_G(T)$,
hence~$N_G(T)=O(N_H(T))\leftthreetimes S$.
Since $N_G(S)\leq N_G(T)$, we obtain that the set of elements of odd order
is a normal subgroup of $N_G(S)$, i.~e., every element of odd order of
$N_G(S)$ is contained in $O(N_G(S))$. On the other hand $S$ is normal in
$N_G(S)$ by definition and $S\cap O(N_G(S))=\{1\}$, hence
$N_G(S)=S\times O(N_G(S))=SC_G(S)$.
\end{proof}

\section{Groups of Lie type}

Our notations for groups of Lie type agrees with \cite{Car1} and for
linear algebraic groups agrees with~\cite{Hu1}. If $G$ is a finite
group of Lie type with trivial centre (we do not exclude
non-simple groups of Lie type, such as $A_1(2)$, all exceptions are given in
\cite[Theorems 11.1.2 and~14.4.1]{Car1}), then $\widehat{G}$
denotes the group of inner-diagonal automorphisms of $G$. In view of
\cite[3.2]{Ste2} we have that $\Aut(G)$ is generated by
inner-diagonal, field and graph automorphisms. Since we are
assuming that $Z(G)$ is trivial, we have that $G$ is isomorphic to the
group of its inner automorphisms and hence we may suppose
that~$G\leq\widehat{G}\leq \Aut(G)$.

Let $\overline{G}$ be a simple connected linear algebraic group over
an algebraically closed field $\mathbb{F}_p$ of positive
characteristic $p$. It is possible here that $Z(\ov{G})$ is nontrivial. An
automorphism $\sigma$ of $\ov{G}$ is called a
{\em Frobenius map} if
$\ov{G}_\sigma$ is finite. Groups $O^{p'}(\ov{G}_\sigma)$ are called {\em
canonical
finite groups of Lie type} and every group $G$ satisfying
$O^{p'}(\ov{G}_\sigma)\leq G\leq \ov{G}_\sigma$ is called a {\em
finite group of Lie type}. Note that in \cite{Car1} only groups
$O^{p'}(\ov{G})$ are called groups of Lie type. But later in \cite{Car6}
R.Carter said that every group $\ov{G}_\sigma$ is a finite group of Lie type
for an arbitrary connected reductive group $\ov{G}$. More over, in \cite{Ca2}
and
\cite{der1} every group $G$ with $O^{p'}(\ov{G}_\sigma)\leq G\leq
\ov{G}_\sigma$ is called a finite group of Lie type. Thus, by given definition
of finite groups of Lie type and canonical finite groups of Lie type we intend
to clarify the situation here. For example, $\P SL_2(3)$ is a canonical finite
 group of Lie type and $\P GL_2(3)$ is a finite group of Lie type. Note that an
element of order $3$ is not conjugate to its inverse in $\P SL_2(3)$ and is
conjugate to its inverse in $\P GL_2(3)$. Since such information about conjugation
is important in many cases (and is very important and useful in this paper), we find
it reasonable to use such notation.

We say that groups
${}^2A_n(q^2)$,
${}^2D_n(q^2)$, ${}^2E_6(q^2)$ are defined over $GF(q^2)$, groups
${}^3 D_4(q^3)$ are defined over $GF(q^3)$ and over groups are defined
over $GF(q)$. The field $GF(q)$ in all cases is called the {\em base
field}. In view of \cite[12.3]{Ste} and \cite[Exercise after Lemma 58]{Ste1} we
have that if $\ov{G}$
is of adjoint type then $\ov{G}_\sigma$ is a group of inner-diagonal
automorphisms of $O^{p'}(\ov{G}_\sigma)$.  If $\ov{G}$ is simply
connected, then $\ov{G}_\sigma=O^{p'}(\ov{G}_\sigma)$
(cf. \cite[12.4]{Ste}). In general for given finite group of
Lie type $G$ (if we consider it as an abstract group) the
corresponding algebraic group is not uniquely determined. For example, if $G=
\P SL_2(5)\simeq SL_2(4)$, then $G$ can be obtained either as
$(SL_2(\F_2))_\sigma$, or as $O^{5'}((\P SL_2(\F_5))_\sigma)$ (for appropriate
$\sigma$). So, for any
finite group of Lie type $G$, we fix (in some way) corresponding algebraic
group $\ov{G}$ and a Frobenius map~$\sigma$ such that
$O^{p'}(\ov{G}_\sigma)\leq G\leq \ov{G}_\sigma$. Let $U=\langle X_r\vert r\in
\Phi^+\rangle$ be the maximal unipotent
subgroup of $G$. If we fix an order on $\Phi(\ov{G})$ consistent with the sum
of roots, then every $u\in U$ can be uniquely written as
\begin{equation}\label{canonicalform}
u=\prod_{r\in\Phi^+}x_r(t_r),                           \end{equation}
where roots are taken in given order and $t_r$ are from the field of definition
of~$G$. Sometimes we use notation $\Phi^\varepsilon(q)$, where
$\varepsilon\in\{+,-\}$, and $\Phi^+(q)=\Phi(q)$ is a split group of Lie type
with base field $GF(q)$, $\Phi^-(q)={}^2\Phi(q^2)$ is a twisted group of Lie
type defined over a field $GF(q^2)$ (with base field~$GF(q)$).

Now let $\ov{R}$ be a closed $\sigma$-stable subgroup of $\ov{G}$. Then we can
consider $R=G\cap \ov{R}$ and $N(G,R)=G\cap N_{\ov{G}}(\ov{R})$. Note that
$N(G,R)\not=N_G(R)$ in general and we call $N(G,R)$ an {\em algebraic normalizer} of
$R$. For example, if we consider $G=SL_n(2)$, then the group of diagonal matrices
$H$ of $G$ is trivial, hence $N_G (H)=G$. But $G=(SL_n(\F_2))_\sigma$, where
$\sigma$ is the Frobenius map $\sigma:(a_{i,j})\mapsto (a_{i,j}^2)$. Then
 $H=\ov{H}_\sigma$, where $\ov{H}$ is the subgroup of diagonal matrices in
$SL_n(\F_2)$. Thus $N(G,H)$ is the group of monomial matrices of $G$. We use term
``algebraic normalizer'' in order to avoid such difficulties,  and to make our
proofs universal. A group $R$ is said to be a {\em torus} (resp. a {\em reductive
subgroup}, a {\em parabolic subgroup}, a {\em maximal torus}, a {\em reductive
subgroup of maximal rank}) if $\ov{R}$ is a torus (resp. a reductive subgroup, a
parabolic subgroup, a maximal torus, a reductive subgroup of maximal rank) of
$\ov{G}$. If $\ov{R}$ is a connected reductive subgroup of maximal rank of $\ov{G}$,
then $\ov{R}=\ov{G}_1\ast\ldots\ast \ov{G}_k\ast\ov{S}$, where $\ov{G}_i$ is a
simple connected linear algebraic group and $\ov{S}=Z(\ov{R})^0$ (see
\cite[Theorem~27.5]{Hu1}). Moreover, if $\Phi_1,\ldots,\Phi_k$ are root systems of
$\ov{G}_1,\ldots,\ov{G}_k$ respectively, then $\Phi_1\oplus\ldots\oplus\Phi_k$ is a
subsystem of $\Phi(\ov{G})$. There is a nice algorithm due to Borel and de Siebental
\cite{BorSie} and independently Dynkin \cite{Dyn} of determining subsystems of
$\Phi$. One has to remove some nodes from the extended Dynkin diagram of $\Phi$. The
remaining connected components are Dynkin diagrams of indecomposable components in
some subsystem and any subsystem can be derived in this way.

Now assume that $\ov{R}$ is $\sigma$-stable. In view of
\cite[10.10]{Ste} there exists a $\sigma$-stable maximal torus
$\ov{T}$ of $\ov{R}$. Let $\ov{G}_1,\ldots,\ov{G}_l$ be the
$\sigma$-orbit of $\ov{G}_1$. Then
$$(\ov{G}_1\ast\ldots\ast\ov{G}_l)_\sigma=
\{x\in\ov{G}_1\vert x= g\cdot g^\sigma\cdot\ldots\cdot
g^{\sigma^{l-1}}\text{ for some }g\in\ov{G}\}_\sigma\simeq
(\ov{G}_1)_{\sigma^l}.$$ In view
of \cite[10.15]{Ste} we have that $\ov{G}_{\sigma^l}$ is finite, hence
$O^{p'}((\ov{G}_1)_{\sigma^l})$ is a canonical finite group of Lie type,
probably, with the base field larger then the base field of
$O^{p'}(\ov{G}_\sigma)$. Since $\ov{G}_1\ast\ldots\ast\ov{G}_l$ is
$\sigma$-stable then $\ov{G}_1\ast\ldots\ast\ov{G}_l\cap \ov{T}$ is a
$\sigma$-stable maximal torus of $\ov{G}_1\ast\ldots\ast\ov{G}_l$.  Therefore we
may assume
that  for
any $\sigma$-orbit $\{\ov{G}_{j_1},\ldots,\ov{G}_{j_i}\}$
$\ov{T}\cap \ov{G}_{j_1}\ast\ldots\ast\ov{G}_{j_i}$ is a maximal $\sigma$-stable
torus of $\ov{G}_{j_1}\ast\ldots\ast\ov{G}_{j_i}$. For $i^{\text{th}}$
$\sigma$-orbit $\ov{G}_{j_1}\ast\ldots\ast\ov{G}_{j_i}$ let
$G_i=O^{p'}(\ov{G}_{j_1}\ast\ldots\ast\ov{G}_{j_i})_\sigma
=O^{p'}(G_{j_1})_{\sigma^i}$. Then
$\ov{R}_\sigma=\ov{T}_\sigma(G_1\ast\ldots\ast G_m\ast
\ov{S}_\sigma)$ and $\ov{T}_\sigma$ normalizes each of $G_i$. Subgroups $G_i$ of
$O^{p'}(\ov{G}_\sigma)$ arising in
this way we call {\em subsystem subgroups} of $O^{p'}(\ov{G}_\sigma)$.

For a $\sigma$-orbit $\{\ov{G}_{j_1},\ldots,\ov{G}_{j_i}\}$ of
$\ov{G}_{j_1}$, with $G_i=O^{p'}((\ov{G}_{j_1})_{\sigma^s})$,
consider $\Aut_{\ov{R}_\sigma}(G_i)$. Since
$G_1\ast\ldots\ast G_{i-1}\ast G_{i+1}\ast\ldots\ast
G_k\ast\ov{S}_\sigma\leq C_{\ov{R}_\sigma}(G_i)$, we have that
$\Aut_{\ov{R}_\sigma}(G_i)\simeq
\left(\left(\ov{T}\ov{G}_{j_1}\right)/
Z\left(\ov{T}\ov{G}_{j_1}\right)
\right)_{\sigma^i},$ i.~e.,
$\Aut_{\ov{R}_\sigma}(G_i)$ is a finite
group of Lie type and $\Aut_{\ov{R}_\sigma}(G_i)$ has trivial
centre. Therefore we may assume that $\P G_i\leq
\Aut_{\ov{R}_\sigma}(G_i)\leq \widehat{\P G_i}.$

Let $\ov{R}$ be a $\sigma$-stable connected reductive  subgroup of maximal rank
(in
particular, $\ov{R}$ can be a maximal torus) of $G$. Let
$Cl(\ov{G}_\sigma,\ov{R})$ be the set of $\ov{G}_\sigma$-conjugacy classes of
$\sigma$-stable subgroups $\ov{R}^g$, where $g\in\ov{G}$. Then
$Cl(\ov{G}_\sigma, \ov{R})$ is in $1-1$ correspondence with the set of
$\sigma$-conjugacy classes $Cl(N_W(W_{\ov{R}})/W_{\ov{R}},\sigma)$ (we define
this term below), where $W$
is the Weyl group of $\ov{G}$, $W_{\ov{R}}$ is the Weyl group of $\ov{R}$ (and
it is a subgroup of $W$). Now define $Cl(N_W(W_{\ov{R}})/W_{\ov{R}},\sigma)$.
Since $N_{\ov{G}}(\ov{R})/\ov{R}\simeq
N_W(W_{\ov{R}})/W_{\ov{R}}$ we obtain the induced action of $\sigma$ on
$N_W(W_{\ov{R}})/W_{\ov{R}}$ and we say that $w_1\equiv w_2$, for $w_1,w_2\in
N_W(W_{\ov{R}})/W_{\ov{R}}$  if there exists
$w\in W$ with $w_1=w^{-1}w_2w^\sigma$. Now if $w$ is an element of
$N_W(W_{\ov{R}})/W_{\ov{R}}$, and $(\ov{R}^g)_\sigma$ corresponds to
the $\sigma$-conjugated class of $w$  then we say that $(\ov{R}^g)_\sigma$ is
obtained by ``twisting'' $\ov{R}$ with $w\sigma$. For more details
see~\cite{Ca5}.

\begin{lem}\label{CentrOfInvolution}
Let $G$ be a simple connected linear algebraic group over a field of
characteristic $p$. Let $t$ be an element of prime order $r\not=p$ of $G$.

Then $C_G(t)/(C_G(t)^0)$ is an $r$-group.
\end{lem}

\begin{proof}
Since $r$ is distinct from the characteristic it follows that $t$ is semisimple.
Hence, $C_G(t)^0$ is a connected reductive subgroup of maximal rank of $G$ and every
$p$-element of $C_G(t)$ is contained in $C_G(t)^0$. Assume that some prime $s\not=r$
divides $\vert C_G(t)/(C_G(t)^0)\vert$. Then $s\not=p$ and $C_G(t)$ contains an
element $x$ of order $s^k$ such that $x\not\in C_G(t)^0$. Since $x,t$ commute we
have that $x\cdot t$ is a semisimple element of $G$. Therefore there exists a
maximal torus $T$ of $G$ with $x\cdot t\in T$. Then $(xt)^r=x^r\in T$. Since
$(s,r)=1$  we have that there exists $m$ such that $rm\equiv 1\pmod{s^k}$, thus
$(x^r)^m=x\in T$. But $T\leq C_G(t)^0$, hence $x\in C_G(t)^0$, a contradiction.
\end{proof}

Assume now that $\ov{R}$ is a $\sigma$-stable parabolic subgroup of
$\ov{G}$. Then it has the unipotent radical $\ov{U}$ and a connected
reductive subgroup $\ov{L}$ such that $\ov{R}/\ov{U}\simeq \ov{L}$. The
subgroup
$\ov{L}$ is called a {\em Levi factor} of $\ov{R}$. Moreover, if
$\ov{S}=Z(\ov{L})^0$, then $\ov{L}=C_{\ov{G}}(\ov{S})$. Let
$Rad(\ov{R})$ be the radical of $\ov{R}$. Then it is a $\sigma$-stable
connected solvable subgroup, hence, by \cite[10.10]{Ste} it contains a
$\sigma$-stable torus $\ov{S}$. Now $C_{\ov{G}}(\ov{S})=C_{\ov{R}}(\ov{S})$ is a
$\sigma$-stable Levi factor of $\ov{R}$, i.~e., every $\sigma$-stable
parabolic subgroup of $\ov{G}$ contains a $\sigma$-stable Levi factor
$\ov{L}$ and $\ov{L}$ is a  connected reductive subgroup of maximal
rank of~$\ov{G}$.

\begin{lem}\label{centUH}
Let $O^{p'}(\overline{G}_\sigma)\leq G\leq\overline{G}_\sigma$ be a finite group of
Lie type over a field of odd characteristic $p$ and the root system $\Phi$ of
$\overline{G}$  be one of the following: $A_\ell$ $(\ell\ge2)$, $D_\ell$ $(\ell\ge
3)$, $B_\ell$ $(\ell\ge 3)$, $E_6$, $E_7$ or $E_8$ and $G\not\simeq{}^3D_4(q^3)$.
Let $U$ be a maximal unipotent subgroup of $G$ and let $H$ be a Cartan subgroup of
$G$ which normalizes $U$.  Then $C_U(\Omega(H))=\{1\}$, where $\Omega(H)=\{h\in H
\mid h^2=1\}$.
\end{lem}

\begin{proof}
For split case and twisted case with $\Phi=D_\ell$, the lemma is proven in
\cite[Lemma~2.8]{TamVdo}. The remaining cases can be proven by using
the same arguments.
\end{proof}

\begin{lem} \label{centUHsymp}
Let $O^{p'}(\overline{G}_\sigma)= G$ be a canonical finite
group of Lie type over a field of odd characteristic $p$ and $-1$ is not a
square in the base field of $G$. Assume that the root system
$\Phi$ of $\overline{G}$  is equal to $C_\ell$.  Let $U$
be a maximal unipotent subgroup of $G$ and let $H$ be a Cartan subgroup of $G$
which normalizes $U$.  Then $C_U(\Omega(H))=\langle X_r\vert r\textrm{ is a long
root}\rangle$, where
$\Omega(H)=\{h\in H \mid h^2=1\}$.
\end{lem}

\begin{proof}
If $r$ is a short root, then there exists a root $s$ with $<s,r>=1$. Thus
$x_{r}(t)^{h_{s}(-1)}=x_{r}((-1)^{<s,r>}t)=x_r(-t)$ (cf.
\cite[Proposition~6.4.1]{Car1}). Therefore, if $x\in C_U(\Omega(H))$ and $x_r(t)$ is
a nontrivial multiplier in decomposition \eqref{canonicalform} of $x$, then $r$ is a
long root. Now if $r$ is a long root, then, for every root $s$, either $\vert
<s,r>\vert=2$, or $<s,r>=0$, i.~e., $x_r(t)^{h_s(-1)}=x_r(t)$. Under our conditions
$\langle h_s(-1)\vert s\in\Phi\rangle=\Omega(H)$, and the lemma follows.
\end{proof}

The following lemma is immediate from~\cite[Theorem~1]{Kond}.

\begin{lem}\label{Norm2Syl}
Let $O^{p'}(\overline{G}_\sigma)= G$ be a canonical
finite group of Lie type and $\overline{G}$ is either of type $A_n$ or  of
type $C_n$, $p$ is odd, $q=p^\alpha$ is the order of the base field of
$G$, and $G$ is split. Let $S$ be a Sylow $2$-subgroup of $G$.

Then  $N_G(S)=SC_G(S)$ if and only if~$q\equiv\pm1\pmod{8}$.
\end{lem}

\begin{lem}\label{InvolutionsAndTori}
Let $O^{p'}(\ov{G}_\sigma)\leq G\leq\ov{G}_\sigma$ be a  finite group of Lie
type with the base field
of characteristic $p$ and order $q$, let
$\ov{G}$ be of adjoint type. Assume also that $G$ is not
isomorphic to ${}^2D_{2n}(q^2)$, ${}^3D_4(q^3)$, ${}^2B_2(2^{2n+1})$,
${}^2G_2(3^{2n+1})$, ${}^2F_4(2^{2n+1}).$ Then there exists a maximal
$\sigma$-stable torus $\ov{T}$ of $\ov{G}$ such that
\begin{itemize}
\item[{\em (1)}] $(N_{\ov{G}}(\ov{T})/\ov{T})_\sigma\simeq
(N_{\ov{G}}(\ov{T}))_\sigma/(\ov{T}_\sigma)=N(\ov{G}_\sigma,\ov{T}_\sigma)/\ov{T
}_\sigma \simeq W$, where $W$ is the Weyl
group of $\ov{G}$;
\item[{\em (2)}] if $r$ is a prime divisor of $q-(\varepsilon1)$, where
$\varepsilon=+$, if $G$ is split and $\varepsilon=-$ if $G$ is twisted, then, up
to conjugation in $O^{p'}(\ov{G}_\sigma)$, every
element of order $r$  is contained in~$\ov{T}_\sigma$;
\item[{\em (3)}] torus $\ov{T}$ is unique, up to conjugation
in~$O^{p'}(\ov{G}_\sigma)$.
\end{itemize}
 \end{lem}

 \begin{proof}
 Since for every maximal torus $T$ of $\ov{G}_\sigma$ we have that
$\ov{G}_\sigma=TO^{p'}(\ov{G}_\sigma)$, without lost we may assume that
$G=\ov{G}_\sigma$.  If $G$ is split then the lemma is evident. In this case $\ov{T}$
is a maximal torus such that $\ov{T}_\sigma$ is a Cartan  subgroup of
$\ov{G}_\sigma$ and (1) is clear. By \cite[F, \S6]{BorCar} we have that every
element of order $r$ of $\ov{G}_\sigma$, up to conjugation, is contained in
$\ov{T}_\sigma$ and (2) follows.  By information about the classes of maximal tori
given in \cite[G]{BorCar} and \cite{Ca3} we have that $\ov{T}$ is unique, up to
conjugation in~$G$.

Assume that $G\simeq {}^2A_n(q^2)$. Then $\ov{T}$ is a maximal torus such that
$\vert \ov{T}_\sigma\vert=(q+1)^n$. The uniqueness follows from
\cite[Proposition~8]{Ca2}. Direct calculations by using
\cite[Proposition~3.3.6]{Car6} show that
$N(\ov{G}_\sigma,\ov{T}_\sigma)/\ov{T}_\sigma\simeq W(\ov{G})=\Sym_{n+1}.$
Assume that $t$ is an element of order $r$ in $G$ (recall that in this case $r$
divides $q+1$). Let $S=A_n(q^2)$ be chosen so that
$G=S_\tau$ for some automorphism $\tau$ (the existence of such group
follows from \cite[Chapter~13]{Car1}). Moreover there exists a Frobenius map
$\rho$ of $\ov{G}$ such that $S=\ov{G}_\rho$ and $\rho=\sigma^2$. Now
let $\ov{H}$ be a $\sigma$-stable maximal torus of $\ov{G}$ such that
$\ov{H}_\sigma$ is a Cartan subgroup of $G$. Then $\ov{H}$ is also
$\rho$-stable and $\ov{H}_\rho$ is a Cartan subgroup of $S$. In view
of \cite[F, \S6]{BorCar} we have that $t$ is contained in $\ov{H}_\rho$. We
state that, up to
conjugation in $S$, the torus $\ov{T}_\sigma$ is contained in
$\ov{H}_\rho$. Indeed,  $\ov{T}_\sigma$ is obtained from $\ov{H}$ by
``twisting'' with $w_0\sigma$, where $w_0\in W(\ov{G})$ is the unique element
that maps
all positive roots onto negative roots. Now $\ov{T}_\rho$ is obtained from
$\ov{H}$ by ``twisting'' with an element $(w_0\sigma)^2=w_0^2\rho=\rho$, i.~e.,
$\ov{T}_\rho$ and $\ov{H}_\rho$ are conjugate in $S$. Let $r_1,\ldots,r_n$ be
the set of fundamental roots of $A_n$. Then $t$, as an element of $\ov{H}_\rho$
can be written as $h_{r_1}(\zeta_1)\cdot\ldots\cdot
h_{r_n}(\zeta_n)$. Now $\ov{T}_\sigma=(\ov{H}_\rho)_{\tau
w_0}$. But $\tau w_0:h_r(\lambda)\mapsto h_{-r}(\lambda^q)=h_r(\lambda^{-q})$,
i.~e. $t^{\tau w_0}=t^{-q}$. Now assume that $t$ is of order $r$. Since $r$
 divides $q+1$ we obtain that $t^{q+1}=e$, i.~e., $t=t^{-q}$. Hence $t^{\tau
w_0}=t$ and~${t\in \ov{T}_\sigma}$.

For $G={}^2D_{2n+1}(q^2)$ we take $\ov{T}$ to be the unique (up to conjugation in
$G$) maximal torus of order $\vert \ov{T}_\sigma\vert= (q+1)^{2n+1}$ (the uniqueness
follows from \cite[Proposition~10]{Ca2}) and for $G={}^2 E_6(q^2)$ we take $\ov{T}$
to be the unique (again up to conjugation in $G$) maximal torus of order $\vert
\ov{T}_\sigma\vert=(q+1)^6$ (the uniqueness follows from \cite[Table~1,
p.~128]{der1}). Like in case $G={}^2 A_n(q^2)$ it is easy to show that $\ov{T}$
satisfies (1) and (2) of the lemma.
 \end{proof}

 \begin{lem}\label{NormOfRegularElementIsNotCentr}
 Let $G$ be a finite group of Lie type and $\ov{G}$, $\sigma$ are chosen so
that $O^{p'}(\ov{G}_\sigma)\leq G\leq \ov{G}_\sigma$. Let $s$ be a regular
semisimple element of odd prime order of $G$.

Then $N_G(C_G(s))\not= C_G(s)$.
 \end{lem}

 \begin{proof}
 In view of \cite[F, \S4 and Proposition~5]{BorCar} we have that
$C_{\ov{G}}(s)/C_{\ov{G}}(s)^0$ is isomorphic to a subgroup of $\Delta$. Now,
if the root system $\Phi$ of $\ov{G}$ is not equal to either $A_n$, or $E_6$,
then $\vert \Delta\vert$ is a power of $2$. Thus, Lemma \ref{CentrOfInvolution}
implies that $C_{\ov{G}}(s)=C_{\ov{G}}(s)^0=\ov{T}$ and
$C_G(s)=C_{\ov{G}}(s)\cap G=T$. Since $N_G(T)\geq N(G,T)\not=T$ we obtain the
statement in this case. Thus we may assume that either $\Phi= A_n$,
or~${\Phi=E_6}$.

Assume first that $\Phi=A_n$, i.~e., $G=A_n^\varepsilon(q)$, where
$\varepsilon\in\{+,-\}$.   Clearly $T=C_{\ov{G}}(s)^0\cap G$ is a normal subgroup of
$C_G(s)$, hence $C_G(s)\leq N(G,T)$.  Assume that $N_G(C_G(s))=C_G(s)$. Then
$C_G(s)=N_{N(G,T)}(C_G(s))$ and $C_G(s)/T$ is a self-normalizing subgroup of
$N(G,T)/T$. As we noted above $C_G(s)/T$ is isomorphic to a subgroup of $\Delta$,
i.~e., it is cyclic. By Lemma \ref{CentrOfInvolution}, we also have that $C_G(s)/T$
is an $r$-group, thus $C_G(s)/T=\langle x\rangle$ for some $x\in N(G,T)/T$ and
$\langle x\rangle$ is a Carter subgroup of $N(G,T)/T$.  Now, in view of
\cite[Proposition~3.3.6]{Car6}, we have that $N(G,T)/T\simeq C_{\Sym_{n+1}}(y)$ for
some $y\in \Sym_{n+1}$. Clearly $C_{C_{\Sym_{n+1}}(y)}(x)$ contains $y$, thus $y$
must be an $r$-element, otherwise $N_{C_{\Sym_{n+1}}(y)}(\langle x\rangle)$ would
contain an element of order coprime to $r$, i.~e., $N_{C_{\Sym_{n+1}}(y)}(\langle
x\rangle)\not=\langle x\rangle$. A contradiction with the fact that $\langle
x\rangle$ is a Carter subgroup of $C_{\Sym_{n+1}}(y)$.

Now let $y=\tau_1\cdot\ldots\cdot \tau_k$ be the decomposition of $y$ into the
product of independent cycles and $l_1,\ldots,l_k$ be the lengths of
$\tau_1,\ldots,\tau_k$ respectively. Assume that first $m_1$ cycles has the same
length $l_1$,  $m_2$ cycles has the length $l_2$ etc. Let
$m_0=n+1-(l_1+\ldots+l_k)$. Then $$C_{\Sym_{n+1}}(y)\simeq
\left(\left(Z_{l_1}\times\ldots\times Z_{l_k}\right)\leftthreetimes
\left(\Sym_{m_1}\times\Sym_{m_2}\times\ldots\right)\right)\times \Sym_{m_0},$$ where
$Z_{l_i}$ is a cyclic group of order $l_i$. If $m_j>1$ for some $j\ge0$, then there
exists a normal subgroup $H$ of $C_{\Sym_{n+1}}(y)$ such that
$C_{\Sym_{n+1}}(y)/H\simeq \Sym_{m_j}\not=\{e\}$. In view of \cite[Table]{Vdo} and
\cite[Table]{TamVdo} we obtain that Carter subgroup in group $S$ satisfying
$\mathrm{Alt}_\ell\leq S\leq \Sym_\ell$ are conjugate for all $\ell\ge1$. Thus
$C_{\Sym_{n+1}}(y)$ and $H$ satisfy ($*$) and $\langle x\rangle$ is the unique, up
to conjugation, Carter subgroup of $C_{\Sym_{n+1}}(y)$. By Lemma
\ref{HomImageOfCarter} we obtain that $\langle x\rangle$ maps onto a Carter subgroup
of $C_{\Sym_{n+1}}(y)/H\simeq \Sym_{m_j}$. In view of \cite{DT1} we have that only a
Sylow $2$-subgroup of $\Sym_{m_j}$ can be a Carter subgroup of $\Sym_{m_j}$. A
contradiction with the fact that $x$ is an $r$-element and $r$ is odd.

Thus we may assume that $C_{\Sym_{n+1}}(y)=\left(Z_{l_1}\times\ldots\times
Z_{l_k}\right)$ and $l_i\not=l_j$ if $i\not=j$. From the known structure of
maximal tori and their normalizers of $A_n^\varepsilon(q)$ (cf.
\cite[Propositions~7,8]{Ca2}, for example) we obtain that
$T=((T_1\times\ldots\times T_k)/Z)\cap A_n^\varepsilon(q)$, where $T_i$ is a,
so-called, Singer group of $GL_{l_i}^\varepsilon(q)=G_i$ and
$N(G,T)=((N(G_1,T_1)\times\ldots\times N(G_k,T_k))/Z)\cap A_n^\varepsilon(q)$.
Thus we may assume that  $N(G,T)=C_G(s)$ and $T$ is a
Singer group, i.~e., it is a cyclic group of order $\frac{q^{n+1}-(\varepsilon
1)^{n+1}}{q-(\varepsilon1)}$. It is known that
if $x$ generates $N(G,T)$ modulo $T$, then $T_x$  is in the centre of $G$, therefore,
$x\not\in C_G(s)$ (for details see \cite[p.~187]{Hup}).

In the remaining case $\Phi=E_6$ by direct calculation, using \cite{GAP}, for example, it
is easy to check that, for every $y\in W(E_6)$, $C_{W(E_6)}(y)$ does not contain Carter
subgroups of order $3$. Since $\vert C_G(s)/T\vert$ divides $3$ and $C_G(s)/T$ is a
Carter subgroup in $C_{W(E_6)}(y)$ for some $y$, we obtain a contradiction.
 \end{proof}

\section{Semilinear groups of Lie type}

Now we define some special overgroups of finite groups of Lie
type. First we give precise description of a Frobenius map $\sigma$.
Let $\overline{G}$ be a simple connected linear algebraic group over
an algebraically closed field $\mathbb{F}_p$ of positive
characteristic $p$. Choose a Borel subgroup $\overline{B}$ of
$\overline{G}$, let $\overline{U}=R_u(\overline{B})$ be the unipotent
radical of $\overline{B}$. There exists a Borel subgroup
$\overline{B}^-$ with $\overline{B}\cap\overline{B}^-=\overline{T}$,
where $\overline{T}$ is a maximal torus of $\overline{B}$ (hence of
$\overline{G}$). Let $\Phi$ be the root system of $\overline{G}$ and let
$\{X_r\vert r\in\Phi^+\}$ be the set of $\overline{T}$-invariant
1-dimensional root subgroups of $\overline{U}$. Every $X_r$ is
isomorphic to the additive group of $\mathbb{F}_p$, so every element of
$X_r$ can be written as $x_r(t)$, where $t$ is the image of $x_r(t)$
under above mentioned isomorphism.  Denote by
$\overline{U}^-=R_u(\overline{B}^-)$ the unipotent radical of
$\overline{B}^-$. As above define $\overline{T}$-invariant
1-dimensional subgroups $\{\overline{X}_r\vert r\in\Phi^-\}$ of
$\overline{U}^-$. Then $\overline{G}=\langle
\overline{U},\overline{U}^-\rangle$. Let $\bar{\varphi}$ be a field
automorphism of $\overline{G}$ and $\bar{\gamma}$ be a graph
automorphism of $\overline{G}$. It is known that $\bar{\varphi}$ can
be chosen so that it  acts
by $x_r(t)^{\bar{\varphi}}=x_r(t^p)$ (see \cite[12.2]{Car1} and \cite[1.7]{Car6}, for
example). In view of
\cite[Proposition~12.2.3~and~Proposition~12.3.3]{Car1} we can choose
$\bar{\gamma}$ so that it
acts by $x_r(t)^{\bar{\gamma}}=x_{\bar{r}}(t)$ if $\Phi$ has no roots
of distinct length or by
$x_r(t)^{\bar{\gamma}}=x_{\bar{r}}(t^\lambda)$ for appropriate
$\lambda$ if $\Phi$ has roots of distinct length. Here $\bar{r}$ is
the image of $r$ under the symmetry $\rho$, corresponding to~$\bar{\gamma}$, of
root system $\Phi$. In both cases we can write
$x_r(t)^{\bar{\gamma}}=x_{\bar{r}}(t^\lambda)$, where $\lambda$ is a
field automorphism (probably trivial) of $\mathbb{F}_p$. From these
formulae it is evident that $\bar{\varphi}\cdot
\bar{\gamma}=\bar{\gamma}\cdot\bar{\varphi}$. Let
$n_r(t)=x_r(t)x_{-r}(-t^{-1})x_r(t)$ and $\overline{N}=\langle
n_r(t)\vert r\in\Phi, t\in\mathbb{F}_p\rangle$.  Let
$h_r(t)=n_r(t)n_r(-1)$ and $\overline{H}=\langle h_r(t)\vert r\in\Phi,
t\in\mathbb{F}\rangle$. In view of \cite[Chapters~6 and~7]{Car1},
$\overline{H}$ is a maximal torus of $\overline{G}$,
$\overline{N}=N_{\overline{G}}(\overline{H})$ and $\overline{X}_r$ are
root subgroups with respect to $\overline{H}$. So we can substitute
$\overline{T}$ by $\overline{H}$ and suppose that under our choice
$\overline{T}$ is $\bar{\varphi}$- and $\bar{\gamma}$-
invariant. Moreover $\bar{\varphi}$ induces the trivial automorphism
of~$\overline{N}/\overline{T}$.

Automorphism $\bar{\varphi}^k, k\in\mathbb{N}$ is called a {\em
classical Frobenius automorphism}. We shall call an automorphism
$\sigma$ a {\em Frobenius automorphism} if
$\sigma$ is conjugate under $\ov{G}$ to $\bar{\gamma}^e\bar{\varphi}^k,
e\in\{0,1\}, k\in
\mathbb{N}$. It follows from Lang-Steinberg theorem \cite[Theorem~10.1]{Ste}
that for any
 $\bar{g}\in\ov{G}$, elements $\sigma$ and $\sigma\bar{g}$ are conjugate under
$\ov{G}$. Thus, in view of  \cite[11.6]{Ste}, we have that a Frobenius map, defined
in previous section, coincides with a Frobenius automorphism defined here.

Now fix $\overline{G}$, $\bar{\varphi}$, $\bar{\gamma}$, and
$\sigma=\bar{\gamma}^e\bar{\varphi}^k$; and assume that $\vert\bar\gamma\vert\le 2$,
i.~e., we do not consider the triality automorphism of $D_4$. Consider
$B=\overline{B}_\sigma$, $T=\overline{T}_\sigma$, and $U=\overline{U}_\sigma$. Since
$\overline{B}, \overline{T},$ and $\overline{U}$ are $\bar{\varphi}$- and
$\bar{\gamma}$- invariant, they give us Borel subgroup, Cartan subgroup, and maximal
unipotent subgroup (Sylow $p$-subgroup) of $\overline{G}_\sigma$ (see
\cite[1.7--1.9]{Car6} for details).

Assume first that $e=0$, i.~e., $O^{p'}(\overline{G}_\sigma)$ is not
twisted (is split). Then $U=\langle X_r\vert r\in\Phi^+\rangle$,
where $X_r$ is isomorphic to the additive group of $GF(p^k)=GF(q)$
and every element of $X_r$ can be written in the form $x_r(t), t\in
GF(q)$. Consider also $U^-=\overline{U}^-_\sigma$. As for $U$ we can
write $U^-=\langle X_r\vert r\in \Phi^-\rangle$ and every element of
$X_r$ can be written in the form $x_r(t),t\in GF(q)$. Now we can
define an automorphism $\varphi$ by the restriction of
$\bar{\varphi}$ on $\overline{G}_\sigma$ and automorphism $\gamma$ by
the restriction of $\bar{\gamma}$ on $\overline{G}_\sigma$. By
definition we have that $x_r(t)^\varphi=x_r(t^p)$ and
$x_r(t)^\gamma=x_{\bar{r}}(t^\lambda)$ (see the definition of
$\bar{\gamma}$ above) for all $r\in\Phi$.   Define
$\zeta=\gamma^\varepsilon\varphi^\ell,$ $\varphi^\ell\not=e$,
$\varepsilon\in\{0,1\}$ to be
an automorphism of $\overline{G}_\sigma$ and define
$\bar\zeta=\bar\gamma^\varepsilon\cdot\bar\varphi^\ell$ to be an
automorphism of $\overline{G}$. Choose a $\zeta$-invariant subgroup
$G$ with $O^{p'}(\overline{G}_\sigma)\leq G\leq\overline{G}_\sigma$.
Note that if the root system $\Phi$ of $\ov{G}$ is not $D_{2\ell}$, then
$\ov{G}_\sigma/(O^{p'}(\ov{G}_\sigma))$ is cyclic. Thus for most groups and
automorphisms, except
groups of type $D_{2\ell}$ over a field of odd characteristic,
any subgroup $G$ of
$\overline{G}_\sigma$ satisfying $O^{p'}(\overline{G}_\sigma)\leq
G\leq \overline{G}_\sigma$ is $\gamma$- and $\varphi$-  invariant.
Denote  $\Gamma G=G\leftthreetimes
\langle\zeta\rangle$ and  $\Gamma
\overline{G}=\overline{G}\leftthreetimes \langle\bar\zeta\rangle$.

Assume now that $e=1$, i.~e., $O^{p'}(\overline{G}_\sigma)$ is twisted. Then
$U=\overline{U}_\sigma$ is generated by groups $X_R$, where $$X_R=\langle
\overline{X}_r\vert r\in \{\alpha s+\beta s^\rho\vert\alpha,\beta\ge0,
\alpha,\beta\in\mathbb{Z}\}\cap\Phi^+ \text{ for some }s\in\Phi^+\rangle_\sigma$$
and $\rho$ is the symmetry of Dynkin diagram corresponding to $\bar{\gamma}$,
$U^-=\overline{U}^-$ is generated by groups
$$X_R=\langle \overline{X}_r\vert r\in \{\alpha s+\beta s^\rho\vert\alpha,\beta\ge0,
\alpha,\beta\in\mathbb{Z}\}\cap\Phi^- \text{ for some }s\in\Phi^-\rangle_\sigma.$$
Define $\varphi$ on $U^\pm$ to be the restriction of $\bar{\varphi}$ on $U^\pm$.
Since $O^{p'}(\overline{G}_\sigma)=\langle U^+,U^-\rangle$ we obtain the
automorphism $\varphi$ of $O^{p'}(\overline{G}_\sigma)$. Consider
$\zeta=\varphi^\ell\not=e$ and let $G$ be a $\zeta$-invariant group with
$O^{p'}(\overline{G}_\sigma)\leq G\leq \overline{G}_\sigma$. Then
$\bar{\zeta}=\bar{\varphi}^\ell$ is an automorphism of $\overline{G}$. Define
$\Gamma G=G\leftthreetimes \langle \zeta\rangle$
and~$\Gamma\overline{G}\leftthreetimes\langle \bar{\zeta}\rangle$.

A group $\Gamma G$ defined above is called a {\em semilinear finite group of Lie
type}  (it is called a {\em semilinear canonical finite group of Lie type} if
$G=O^{p'}(\ov{G}_\sigma)$) and $\Gamma \overline{G}$ is called a {\em semilinear
algebraic group}.  Note that $\Gamma\ov{G}$ can not be defined without $\Gamma G$,
since we need to know that $\varphi^\ell\not=e$. If $G$ is written in notations of
\cite{Car1}, i.~e. $G=A_n(q)$ or $G={^2A_n(q^2)}$ etc., then we shall write $\Gamma
G$ by $\Gamma A_n(q)$, $\Gamma {}^2A_n(q^2)$, etc.

Consider $x\in\Gamma G\setminus G$. Then $x=\zeta^k g$ for some
$k\in\mathbb{N}$ and $g\in G$. Define $\bar{x}$ to be
$\bar{\zeta}^kg$. Conversely, if $\bar{x}=\bar{\zeta}^k g$ for some
$g\in G$ and $\zeta^k\not=e$, define $x$ to be equal to~$\zeta^k g$. Note that
we do not
need to suppose that $\bar{x}\in\Gamma \overline{G}\setminus
\overline{G}$ since~$\vert\bar{\zeta}\vert=\infty$. If $x\in G$ we
define~$\bar{x}=x$.

\begin{lem}\label{equivnormalizer}
Let $H$ be a subgroup of $G$. Then $x$ normalizes $H$ if and only if
$\bar{x}$ normalizes $H$ as a subgroup of~$\overline{G}$.
\end{lem}

\begin{proof}
Since $\zeta$ is the restriction of $\bar{\zeta}$ on $G$ our
statement is trivial.
\end{proof}

Let $H_1$ be a subgroup of $\Gamma G$. Then $H_1$ is generated by $H=H_1\cap G$ and
an element $x=\zeta^kg$, moreover $H$ is a normal subgroup of $H_1$. In view of
Lemma \ref{equivnormalizer} we can consider $\overline{H}_1=\langle
\bar{x},H\rangle$. Now we find it reasonable to explain, why we use such complicate
notations and definitions. We have that $\zeta$ is always of finite order, but
$\bar\zeta$ is always of infinite order. Thus, even if $Z(G)$ is trivial, we can not
consider $G\leftthreetimes\langle\bar\zeta\rangle$ as a subgroup of $\Aut(G)$.
Hence, we need to define in some way (one possible way is just given) the connection
between elements of $\Aut(G)$ and elements of $\Aut(\ov{G})$ in order to use the
machinery of linear algebraic groups.

Let $\overline{R}$ be a $\sigma$-stable maximal torus (resp.
reductive subgroup of maximal rank, parabolic subgroup) of
$\overline{G}$, let $y\in N_{\Gamma
\overline{G}}(\overline{R})$ be chosen so that there exists $x\in
\Gamma G$ with $y=\bar{x}$. Then $R_1=\langle
x,\overline{R}\cap G\rangle$ is called a {\em maximal torus} (resp.
a {\em reductive subgroup of maximal rank}, a {\em parabolic
subgroup}) of~$\Gamma G$.

\begin{lem}\label{parabolic}
Let $M=\langle x, X\rangle$ be a subgroup of $\Gamma G$ such that
$X=M\cap G\unlhd M$ and $O_p(X)$ is nontrivial. Then there exists
$\sigma$- and $\bar{x}$- stable parabolic subgroup $\overline{P}$ of
$\overline{G}$ such that $X\leq \overline{P}$ and~$O_p(X)\leq
R_u(\overline{P})$.
\end{lem}

\begin{proof}
Define $U_0=O_p(X)$, $N_0=N_{\overline{G}}(U_0)$. Then
$U_i=R_u(N_{i-1})$ and $N_i=N_{\overline{G}}(U_i)$. Clearly $U_i$,
$N_i$ are $\bar{x}$- and $\sigma$- stable for all $i$. In view of
\cite[Proposition~30.3]{Hu1}, the chain of subgroups $N_0\leq N_1\leq
\ldots\leq N_k\leq\ldots$ is finite and $\overline{P}=\cup_i N_i$ is
a proper parabolic subgroup. Clearly $\overline{P}$ is $\sigma$- and
$\bar{x}$- stable.
\end{proof}

\begin{lem}\label{Syl2InCentrOfFieldAut}
Let $G$ be a finite  group of Lie type over a field of odd
characteristic $p$ and $G\not\simeq
{}^2G_2(3^{2n+1})$, ${}^3D_4(q^3)$, ${}^2D_{2n}(q^2)$. Assume that $\ov{G}$
and $\sigma$ are chosen so that $O^{p'}(\overline{G}_\sigma)\leq
G\leq\overline{G}_\sigma$. Let $\psi$ be a field
automorphism of $O^{p'}(\overline{G}_\sigma)$ of odd order.

Then a Sylow $2$-subgroup of $G_\psi$ is a Sylow $2$-subgroup of~$G$. Moreover
there exists a maximal torus $T$ of $G$ such that $N(G,T)/T\simeq
N_{\ov{G}}(\ov{T})/\ov{T}$, a Sylow $2$-subgroup of $T_\psi$ is a Sylow
$2$-subgroup of~$T$, and $\bar\psi$ normalizes every $\ov{T}$-invariant root
subgroup $\ov{X}$ of~$\ov{G}$.
\end{lem}

\begin{proof}
 Assume that
$\vert \psi\vert=k$. Let $GF(q)$ be the base field of $G$. Then $q=p^\alpha$ and
$\alpha=k\cdot m$. It is easy to check that every field automorphism of odd order
centralizes the Sylow $2$-subgroup of $\ov{G}_\sigma/(O^{p'}(\ov{G}_\sigma))$, hence
we may assume that $G=\ov{G}_\sigma$. Now $\vert G\vert$ can be written as $\vert
G\vert=q^N(q^{m_1}+\varepsilon_11)\cdot\ldots\cdot(q^{m_n} +\varepsilon_n 1)$ for
some $N$, where $n$ is the rank of $G$, $\varepsilon_i=\pm$.  (cf.
\cite[Theorems~9.4.10 and~14.3.1]{Car1}). Similarly we have that $\vert
G_\psi\vert=(p^m)^N((p^m)^{m_1}+\varepsilon_11)\cdot\ldots\cdot
((p^m)^{m_n}+\varepsilon_n1)$, i.~e., $\vert G\vert_2=\vert G_\psi\vert_2$ and a
Sylow $2$-subgroup of $G_\psi$ is a Sylow $2$-subgroup of $G$.

Now, by Lemma \ref{InvolutionsAndTori}, there exists a maximal torus $T$ of $G_\psi$
such that $N(G_\psi,T)/T\simeq N_{\ov{G}}(\ov{T})/\ov{T}$ and $\vert
T_\psi\vert=(p^m-\varepsilon1)^n$. Clearly $\vert\ov{T}\cap
G\vert=(q-\varepsilon1)^n$. If $G$ is split, then $T_\psi$ is a Cartan subgroup of
$G_\psi$ and the lemma is evident. Thus we may assume that $G$ and $G_\psi$ are
twisted, in particular, $\varepsilon=-$. In view of proof of Lemma
\ref{InvolutionsAndTori}, we have that there exists a split group $L$ such that
$G_\psi=L_\lambda$ for some automorphism $\lambda$ of order $2$, $\ov{T}\cap L$ is a
Cartan subgroup of $L$, and $\psi$ can be considered as a field automorphism of $L$.
Therefore for every $\ov{T}$-root subgroup $\ov{X}$ of $\ov{G}$, $\ov{X}\cap L$ is a
root subgroup of $L$, and it is $\psi$-invariant. Hence $\ov{X}$ is
$\bar\psi$-invariant.
\end{proof}

\begin{lem}\label{ConjAutomorphisms}
{\em \cite[(7-2)]{GorLyo}} Let $\ov{G}$ be a connected simple linear algebraic
group over a field of characteristic $p$, $\sigma$ be a Frobenius map of
$\ov{G}$ and $G=\ov{G}_\sigma$ be a finite group of Lie type. Let
$\varphi$ be a field or a graph-field automorphism of $G$  and let $\varphi'$ be an
element of $(G\leftthreetimes\langle\varphi\rangle)\setminus G$ such that
$\vert\varphi'\vert=\vert\varphi\vert$.

Then there exists an element $g\in G$ such that
$\langle\varphi\rangle^g=\langle\varphi'\rangle$. In particular, if
$G/O^{p'}(G)$ is a $2$-group and $\varphi$ is of odd order, then such $g$ can
be chosen in~$O^{p'}(G)$.
\end{lem}

The following lemma is proven for classical groups in~\cite{FeZu}.

\begin{lem}\label{ConjInverseInGraph}
Let $G$ be a  finite group of Lie type with $Z(G)=\{1\}$, $\overline{G}$, $\sigma$
are chosen so that $O^{p'}(\overline{G}_\sigma)\leq G\leq \overline{G}_\sigma$ and
$Z(\overline{G})=\{1\}$. Assume that $\tau$ is  the graph automorphism of order $2$.

Then every semisimple element $s\in G$ is conjugate to its inverse under~$\langle
\tau\rangle\rightthreetimes (O^{p'}(\overline{G}_\sigma))$.
\end{lem}

\begin{proof}
In view of \cite[3.2]{Ste2} we have that if $G$ has a graph automorphism of order
$2$, then $G$ is split. If $\ov{G}$ is not of type $A_n,D_{2n+1},E_6$, then the
lemma follows from \cite[Lemma~2.2]{TamVdo}, thus we need to consider groups of type
$A_n,D_{2n+1},E_6$.  Denote by $\bar\tau$ the graph automorphism of $\ov{G}$ such
that $\bar\tau\vert_G=\tau$. Let $\ov{T}$ be a maximal $\sigma$-stable torus of
$\ov{G}$ such that $\ov{T}_\sigma\cap G$ is a Cartan subgroup of $G$. Let
$r_1,\ldots, r_n$ be fundamental roots of $\Phi(\ov{G})$ and $\rho$ be the symmetry
corresponding to $\bar\tau$. Denote $r_i^\rho$ by $\bar{r}_i$. Then $\ov{T}=\langle
h_{r_i}(t_i)\vert\text{, where } 1\le i \le n\text{ and } t_i\not=0\rangle$ and
$h_{r_i}(t_i)^{\bar\tau}=h_{\bar{r}_i}(t_i)$. Denote by $W$ the Weyl group of
$\ov{G}$. Let $w_0$ be the unique element of $W$ mapping all positive roots onto
negative roots and let $n_0$ be its preimage in $N_{\ov{G}}(\ov{T})$ under the
natural homomorphism $N_{\ov{G}}(\ov{T})\rightarrow N_{\ov{G}}(\ov{T})/\ov{T}\simeq
W$. Since $\sigma$ acts trivially on $W=N(G,T)/T$,  we can take $n_0\in G$, i.~e.,
$n_0^\sigma=n_0$. Then for all $r_i$ and $t$ we have that
$$h_{r_i}(t)^{n_0\bar\tau}=h_{r_i^{w_0\rho}}(t)=h_{-r_i}(t)=h_{r_i}(t^{-1}).$$
Thus $x^{n_0\bar\tau}=x^{-1}$ for all $x\in \ov{T}$.

Now let $s$ be a semisimple element of $G$. Then there exists a maximal
$\sigma$-stable torus $\ov{S}$ of $\ov{G}$ containing $s$. Since all maximal tori of
$\ov{G}$ are conjugate, we have that there exists $g\in \ov{G}$ such that
$\ov{S}^g=\ov{T}$. Therefore $s^{gn_0\bar\tau g^{-1}}=s^{-1}$. Since
$n_0^\sigma=n_0$ and $\bar\tau^\sigma=\bar\tau$ we have that $(gn_0\bar\tau
g^{-1})^\sigma=g^\sigma n_0\bar\tau (g^{-1})^\sigma$. Moreover, since $\ov{S}$ is
$\sigma$-stable, then for every $x\in\ov{S}$ we have that $x^{gn_0\bar\tau
g^{-1}}=x^{g^\sigma n_0\bar\tau (g^{-1})^\sigma}=x^{-1}$, i.~e., $gn_0\bar\tau
g^{-1}\ov{S}=g^\sigma n_0\bar\tau (g^{-1})^\sigma\ov{S}$. In particular, there
exists $t\in\ov{S}$ such that  $gn_0\bar\tau g^{-1}t=g^\sigma n_0\bar\tau
(g^{-1})^\sigma$. In view of Lang-Steinberg Theorem \cite[Theorem~10.1]{Ste} there
exists $y\in\ov{S}$ such that $t=y\cdot (y^{-1})^\sigma$. Therefore, $g n_0\bar\tau
g^{-1}y=(g n_0\bar\tau g^{-1}y)^\sigma$, i.~e., $gn_o\tau g^{-1}y\in
\ov{G}_\sigma\leftthreetimes \langle\tau\rangle$, and $s^{gn_o\tau g^{-1}y}=s^{-1}$.
Since $O^{p'}(\ov{G}_\sigma)\ov{S}_\sigma=\ov{G}_\sigma$, and $\ov{S}_\sigma$ is
Abelian, we may find $z\in \ov{S}_\sigma$ such that $gn_0\tau g^{-1}yz\in
O^{p'}(\ov{G}_\sigma)\leftthreetimes\langle\tau\rangle$.
\end{proof}

\begin{lem}\label{CartBorel}
Let $\Gamma G$ be a finite semilinear group of Lie type over a field of
characteristic $p$ (we do not exclude case $\Gamma G=G$) and $Z(G)=\{e\}$.
Assume that
$B=U\leftthreetimes H$, where $H$ is a Cartan subgroup of $G$, is a
$\zeta$-invariant Borel subgroup of $G$ and that
$B\leftthreetimes \langle\zeta\rangle$ contains a Carter subgroup $K$ of
$\Gamma G$. Consider $K_u=K\cap U$ and assume that $K_u\not=\{e\}$. Then one of
the following holds:
\begin{itemize}
\item[{\em (1)}] Either $\Gamma G={}^2A_2(2^{2t})\leftthreetimes
\langle\zeta\rangle$, or $\Gamma G=\widehat{{}^2A_2(2^{2t})}\leftthreetimes
\langle\zeta\rangle$;  $\vert \zeta\vert=t$ is odd, $C_G(\zeta)\simeq
\widehat{{}^2A_2(2)}$, $K\cap G$ has
order~$2\cdot 3$ and is isomorphic to a Carter subgroup
of~$\widehat{^2A_2(2^2)}$.
\item[{\em (2)}] $G$  is defined over
  $GF(2^t)$, $\zeta$ is a field automorphism,
$\vert\zeta\vert=t$, and $K$ contains a Sylow $2$-subgroup
  of~$G_{\zeta_{2'}}$.
\item[{\em (3)}] $G/Z(G)\simeq \P SL_2(3^{t})$,
  $\vert\zeta\vert=t$ is odd,  and $K$ contains a
  Sylow $3$-subgroup of~$G_{\zeta_{3'}}$.
\item[{\em (4)}] $\Gamma G={}^2G_2(3^{2n+1})\leftthreetimes\langle\zeta\rangle$,
$\vert\zeta\vert=2n+1$, $K\cap {}^2G_2(3^{2n+1})=S\times P$, where $S$ is of
order $2$ and~${\vert P\vert =3^{\vert\zeta\vert_3}}$.
\end{itemize}
\end{lem}

\begin{proof}
If $G$ is one of the following groups: $A_1(q)$, $G_2(q)$, $F _4(q)$,
${}^2B_2(2^{2n+1})$, or ${}^2F_4(2^{2n+1})$, then the lemma follows from
\cite[Table]{Vdo} or \cite[Table]{TamVdo}.  If $\Gamma G=G$ then our result follows
from  \cite{PreTamVdo} and \cite{TamVdo}. So we may assume that $\Gamma G\not=G$,
i.~e., that $\Gamma G$ contains a nontrivial field or a graph-field automorphism
$\zeta$. It is convenient to say that field automorphism of even order of twisted
group is a graph-field automorphism, and through the proof of the lemma we use this
definition. Let $\Phi$ be the root system of~$\ov{G}$. If $\Phi=C_\ell$, then we
prove this lemma later in Theorem~\ref{sympcarter}, so we assume
that~${\Phi\not=C_\ell}$. If $\Phi=D_4$ and either graph-field automorphism contains
a graph automorphism of order $3$, or $G\simeq {}^3D_4(q^3)$ and
$\vert\zeta\vert\equiv0\pmod 3$, then we prove this lemma later in
Theorem~\ref{CarterTriality}. Since we shall use this lemma in the proof of Theorem
\ref{CarterSemilinear}, after Theorems \ref{sympcarter} and \ref{CarterTriality}, it
is possible to make such additional assumptions on~$G$.

Assume that $q$ is odd and $\Phi$ is one of the following types: $A_\ell$
$(\ell\ge2)$, $D_\ell$ $(\ell\ge 3)$, $B_\ell$ $(\ell\ge 3)$, $E_6$, $E_7$ or $E_8$.
By Lemma \ref{HomImageOfCarter} we have that $KU/U$ is a Carter subgroup of
$B\leftthreetimes\langle\zeta\rangle\simeq H\leftthreetimes\langle \zeta\rangle$.
Since $H_\zeta\leq Z(H\leftthreetimes\langle \zeta\rangle)$, we obtain, up to
conjugation in $B$, that $H_\zeta\leq K$. If $\zeta$ is a field automorphism, then
$\Omega(H)\leq H_\zeta$. In view of nilpotency of $K$ we obtain that $K_u\leq
C_U(\Omega(H))$. By Lemma \ref{centUH} it follows that $C_U(\Omega(H))=\{e\}$, a
contradiction with $K_u\not=\{e\}$. If $\zeta$ is a graph-field automorphism, then
$\zeta_2\not=e$ and $\zeta_{2'}$ centralizes $\Omega(H)$. Thus every element of odd
order of $H\leftthreetimes \langle\zeta\rangle$ centralizes $\Omega(H)$ and, up to
conjugation in $B$, we have that $\Omega(H)\leq K$. Again by Lemma \ref{centUH} we
obtain a contradiction.

Assume that $G\simeq{}^2G_2(3^{2n+1})$ and $\Gamma
G=G\leftthreetimes\langle\zeta\rangle$. Again by Lemma \ref{HomImageOfCarter}
we have that $KU/U$ is a Carter subgroup of
$H\leftthreetimes\langle\zeta\rangle$. Since $(2n+1,3^{2n+1}-1)=1$ we have that
$H_\zeta\simeq KU/U\cap HU/U$ is of order $2$. Thus $K\cap
G=K_u\times\langle t\rangle$, where $t$ is an involution. It follows that
$K_u=C_G(t)\cap G_{\zeta_{3'}}$. Now case (4) follows from \cite{war}
and~\cite[Theorem~1]{levnuz}.

Assume now that $q=2^t$ is even. Assume first that $\Phi$ is one of the
following types: $A_\ell$ $(\ell\ge2)$, $D_\ell$ $(\ell\ge 3)$, $B_\ell$
$(\ell\ge 3)$, $E_6$, $E_7$ or $E_8$, $G$ is split, and $\zeta$ is a field
automorphism. It is easy to see that for any $r\in \Phi$ there exists
$s\in\Phi$ such that $<s,r>=1$. If $\vert\zeta\vert\not=t$, then
$h_s(\lambda)\in H_\zeta\leq K$ for some $\lambda\not=1$, hence
$x_r(t)^{h_s(\lambda)}=x_r(t)(\lambda^{<s,r>}t)=x_r(\lambda t)$ (cf.
\cite[Proposition~6.4.1]{Car1}). It follows that $K_u\leq C_U(H_\zeta)=\{e\}$,
a contradiction. Therefore, $\vert\zeta\vert=t$, $H_\zeta=\{e\}$ and we obtain
statement (2) of the lemma.

Now assume that $\Phi$ is  of type: $A_\ell$ $(\ell\ge2)$, $D_\ell$ $(\ell\ge 3)$,
or $E_6$;  and either $G$ is split and $\zeta$ is
a graph-field automorphism, or $ G$ is twisted. Let $\rho$ be the symmetry of the
Dynkin diagram of $\Phi$ corresponding to $\tau$ and denote $r^\rho$ by $\bar{r}$.
If $\Phi\not=A_2$, then for every $r\in \Phi$ there exists $s\in\Phi$ such that
$s+\bar{s}\not\in\Phi$ and $<s,r>=1$. Then we proceed like in case of field
automorphism, taking $h_s(\lambda)h_{\bar{s}}(\lambda^{2^{t/l}})$. If $G=A_2(2^t)$
then $K\geq H_\zeta$ contains an element $x$ of order $3$ such that $x\in
{}^2A_2(2^2)\leq A_2(2^2)\leq A_2(2^t)$. By using \cite{ATLAS} or \cite{GAP} one can
see that $x$ is conjugate to $x^{-1}$ in $A_2(2^2)$, hence in $G$. Since the only
composition factor of $C_G(x)$ is isomorphic to $A_1(2^t)$ (see
\cite[Proposition~7]{Ca2}), then \cite[Table]{Vdo} and \cite[Theorem~3.5]{TamVdo}
imply that $C_{G\leftthreetimes\langle\zeta\rangle}(x)$ satisfies ($*$), a
contradiction with Lemma~\ref{power}.

Now assume that $G\simeq {}^2A_2(2^{2t})$.  By Lemma
\ref{HomImageOfCarter} we have that $KU/U$ is a Carter subgroup of
$H\leftthreetimes\langle\zeta\rangle$.
Now if $\vert\zeta\vert$ is even, then $H_\zeta$  is isomorphic to a Cartan
subgroup of $A_{2}(2^{2t/\vert\zeta\vert})$. If $H_\zeta=\{e\}$, we obtain
statement (2) of the lemma, if $H_\zeta\not=\{e\}$, then $K_u\leq
C_U(H_\zeta)=\{e\}$, and this gives a contradiction with the condition
$N_G(K_u)=B$. If $\vert\zeta\vert\not=t$ is odd, then  $H_\zeta\leq
K$ contains an element $x$ of order greater, than $3$ and direct calculations
show that $C_U(x)=\{e\}$. If $\vert\zeta\vert=t$ is odd, then we obtain statement
(1) of the lemma.
\end{proof}

\section{Carter subgroups in symplectic groups}

From now by $Cmin$ we denote the minimal $n$ such that $A$ is an almost simple
group, $F^\ast(A)$ is a simple group of Lie type of order $n$ and $A$ contains
nonconjugate Carter subgroups. We shall prove that $Cmin=\infty$, i.~e. that
such a group $A$ does not exist.
In this section we consider Carter subgroups
in an almost simple group $A$ with
simple socle $G=F^\ast(A)\simeq \P Sp_{2n}(q)$. We consider such groups here,
since for groups of type $\P Sp_{2n}(q)$ Lemma \ref{centUH} is not true and we
use arguments slightly different from those that we use in proof of
Theorem~\ref{CarterSemilinear}.

\begin{ttt}\label{sympcarter}
Let $G$ be a finite group of Lie type with trivial centre (not
necessary simple) over a field of characteristic $p$  and $\overline{G}$,
$\sigma$ are chosen so that $\P Sp_{2n}(p^t)\simeq
O^{p'}(\overline{G}_\sigma)\leq G\leq\overline{G}_\sigma$.  Choose a subgroup
$A$ of $\Aut(\P Sp_{2n}(p^t))$ containing $G$.  Let $K$ be a Carter subgroup of
$A$. Assume also
that~$\vert \P Sp_{2n}(p^t)\vert\leq Cmin$ and $A=\langle K,G\rangle$.

Then exactly one of the following statements holds:

\begin{itemize}
  \item[{\em (1)}] $G$ is defined over
  $GF(2^t)$, $\vert\zeta\vert=t$, and $K=S\leftthreetimes\langle\zeta\rangle$,
where $S$ is a Sylow $2$-subgroup
  of~$G_{\zeta_{2'}}$.
\item[{\em (2)}] $G\simeq \P SL_2(3^{t})\simeq \P Sp_2(3^t)$,
  $\vert\zeta\vert=t$ is odd,  and $K=S\leftthreetimes\langle\zeta\rangle$,
where $S$ is a
  Sylow $3$-subgroup of~$G_{\zeta_{3'}}$.
\item[{\em (3)}] $p$ does not divide $\vert K\cap G\vert$ and $K$ is
  contained in the
  normalizer of a Sylow $2$-subgroup of~$A$.
\end{itemize}
\end{ttt}

\begin{proof}
Assume by contradiction that $K$ is a Carter subgroup of $\Gamma G$ and $K$ does not
satisfy the theorem. Write $K=\langle x,K_G\rangle$, where~$K_G=K\cap G\unlhd K$. If
either  $p\not= 3$ or $t$ is even, then our result follows from
\cite[Theorem~3.5]{TamVdo}. Thus we may assume that $q=3^t$ and $t$ is odd. By
\cite[Lemma~2.2]{TamVdo} we have that every semisimple element of odd order is
conjugate to its inverse in $G$. Now, for every semisimple element $t\in G$, every
non-Abelian composition factor of $C_G(t)$ is a simple group of Lie type (cf.
\cite{Ca5}) of order less, than $Cmin$. Therefore, for every non-Abelian composition
factor $S$ of $C_G(t)$, Carter subgroups of $\Aut_{C_G(t)}(S)$ are conjugate. It
follows that $C_G(t)$ satisfies $(*)$. Hence, by Lemma \ref{power}, $\vert
K_G\vert=2^\alpha\cdot 3^\beta$ for some $\alpha,\beta\ge 0$.  If $G=\widehat{\P
Sp_{2n}(q)}$, then by \cite[Theorem~2]{Won} we have that every unipotent element is
conjugate to its inverse. Since $3$ is a good prime for $G$, then \cite[Theorem~1.2
and~1.4]{Seitz} imply that, for any element $u\in G$ of order $3$, all composition
factors of $C_G(u)$ are simple groups of Lie type of order less, than $Cmin$. Thus
$C_G(u)$ satisfies ($*$), hence, by Lemma \ref{power}, we obtain that $K_G$ is a
$2$-group. By Lemmas \ref{Syl2InCentrOfFieldAut} and \ref{ConjAutomorphisms} every
element  $x\in A\setminus G$ of odd order with $\langle x\rangle\cap G=\{e\}$
centralizes some Sylow $2$-subgroup of $G$. Hence $K$ contains a Sylow $2$-subgroup
of $A$, i.~e., $K$ satisfies (3) of the theorem. Thus we may assume that $G=\P
Sp_{2n}(q)$ and $\beta\ge 1$, i.~e., a Sylow $3$-subgroup $O_3(K_G)$ of $K_G$ is
nontrivial.  By Lemma \ref{parabolic} we obtain that $K_G$ is contained in some
$K$-invariant parabolic subgroup $P$ of $G$ with the Levi factor $L$ and, up to
conjugation in $P$, a Sylow $2$-subgroup $O_2(K_G)$ of $K_G$ is contained in $L$. We
have that $KO_3(P)/O_3(P)$ is isomorphic to $\widetilde{K}=K/O_3(K_G)$ and, by Lemma
\ref{HomImageOfCarter}, $\widetilde{K}$ is a Carter subgroup of $\langle
\widetilde{K},L\rangle$. Now $\widetilde{K}\cap L=O_2(K_G)$ is a $2$-group and every
element $x\in \langle \widetilde{K},L\rangle\setminus L$ with $\langle x\rangle\cap
L=\{e\}$ of odd order centralizes a Sylow $2$-subgroup of $L$ (cf. Lemmas
\ref{Syl2InCentrOfFieldAut} and \ref{ConjAutomorphisms}). Therefore $O_2(K_G)$
contains a Sylow $2$-subgroup of $L$, in particular, contains a Sylow $2$-subgroup
of $H$. Hence, $K_G$ contains $\Omega(H)$. Since $K$ is nilpotent, Lemma
\ref{centUHsymp} implies that $O_3(K_G)\leq C_U(\Omega(H))=\langle X_r\vert r
\textrm{ is a long root}\rangle$. Since for any two long positive roots $r,s$ we
have that $r+s\not\in\Phi$, Chevalley commutator formulae \cite[Theorem~5.2.2]{Car1}
implies that $\langle X_r\vert r \textrm{ is a long root}\rangle$ is Abelian.

Up to equivalence of root systems, we may suppose that $\Phi$ is contained in a
Euclidean space with orthonormal basis $e_1,\ldots,e_n$, and its roots has the form
$\pm e_i\pm e_j$, $i,j\in\{1,\ldots,n\}$ (short roots) or $\pm 2e_i$,
$i\in\{1,\ldots,n\}$ (long roots). If
$\{r_1,r_2,\ldots,r_{n-1},r_n\}=\{e_1-e_2,e_2-e_3,\ldots,e_{n-1}-e_n,2e_n\}$ is a
set of fundamental roots of $\Phi$, then long positive roots has the following form
$r_n+2r_{n-1}+\ldots+2r_k=2e_k$ for some $k$. Thus there exists a nontrivial
$O_2(K_G)$-invariant subgroup $\langle X_r\vert r\in I\rangle= O_p(P)\cap\langle
X_r\vert r\text{ is a long root}\rangle$, where $I$ is a subset of the set of long
positive roots. Group $O_2(K_G)$ acts by conjugation on $\langle X_r\vert r\in
I\rangle$, thus we obtain a representation $O_2(K_G)\rightarrow \Sym(I)$. Assume
that there exists an orbit $\Omega$ of length greater, than $1$ such that
$O_3(K_G)\cap\langle X_r\vert r\in\Omega\rangle\not=\{e\}$. Without lost we may
assume that it is $\{X_{2e_n},\ldots,X_{2e_k}\}$. Since $K$ is nilpotent, then
$O_3(K_G)\cap \langle X_{2e_n},\ldots,X_{2e_k}\rangle$ contains an element
$v=x_{2e_n}(t)\cdot x_{2e_{n-1}}(t)\cdot\ldots\cdot x_{2e_k}(t)$ for some $t\in
GF(q)$ and it is central in $K$.  Indeed, $K\cap \langle
X_{2e_n},X_{2e_{n-1}},\ldots,X_{2e_k}\rangle$ is normal in $K_G$, $\zeta$ normalizes
$K\cap \langle X_{2e_n},X_{2e_{n-1}},\ldots,X_{2e_k}\rangle$ (since $\zeta$
normalizes each of $X_r$), hence, $K\cap \langle
X_{2e_n},X_{2e_{n-1}},\ldots,X_{2e_k}\rangle$ is normal in $K$. Therefore, $Z(K)\cap
(K\cap \langle X_{2e_n},X_{2e_{n-1}},\ldots,X_{2e_k}\rangle)$ is nontrivial. Since
$O_2(K_G)$ acts transitively on $\Omega$, we obtain required form of $v$. Now,
either $v$, or $v^{-1}$ under $H$ is conjugate to $f=x_{2e_n}(1)\cdot
x_{2e_{n-1}}(1)\cdot\ldots\cdot x_{2e_k}(1)$, therefore we may assume that $v=
x_{2e_n}(1)\cdot x_{2e_{n-1}}(1)\cdot\ldots\cdot x_{2e_k}(1)$. We want to show that
$v$ and $v^{-1}$ are conjugate in $G$. Since $n-k+1$ is even (as the order of an
orbit of a $2$-group), we may write $v=v_k\cdot v_{k+2}\cdot\ldots\cdot v_{n-1}$,
where $v_i=x_{2e_i}(1)x_{2e_{i+1}}(1)$. Now we show that there exist
$x_k,x_{k+2},\ldots,x_{n-1}$ such that
$$v_i^{x_j}=\left\{\begin{array}{ll}v_i^{-1}&\text{ if }i=j,\\ v_i&\text{
if }i\not=j.\end{array}\right. ,$$ i.~e., $v^{x_k\cdot x_{k+2}\cdot\ldots\cdot
x_{n-1}}=v^{-1}$. We construct $x_{n-1}$. We may choose structure constant so
that $C_{1,1,r_{n-1},r_n}=1$, $C_{1,1,r_{n-1},r_{n-1}+r_n}=1$,
$C_{2,1,r_{n-1},r_n}=-1$. Then  Chevalley commutator formulae
\cite[Theorem~5.2.2]{Car1}
implies that
\begin{equation*}
\left(x_{r_n}(1)\cdot
x_{r_n+2r_{n-1}}(1)\right)^{x_{r_{n-1}+r_n}(1)}=
x_{r_n}(1)\cdot x_{r_n+r_{n-1}}(-1)\cdot
x_{r_n+2r_{n-1}}(-1).
\end{equation*}
Consider a reflection $w$ in the root $e_{n-1}-e_n$. We may choose its preimage
$x$ in $G$ so that
\begin{equation*}
(x_{r_n}(1)\cdot x_{r_n+r_{n-1}}(-1)\cdot
x_{r_n+2r_{n-1}}(-1))^x=x_{r_n}(-1)\cdot x_{r_n+r_{n-1}}(-1)\cdot
x_{r_n+2r_{n-1}}(1)
\end{equation*}
At the end, again by using Chevalley commutator formulae, we have
\begin{equation*}
(x_{r_n}(-1)\cdot x_{r_n+r_{n-1}}(-1)\cdot
x_{r_n+2r_{n-1}}(1))^{x_{r_{n-1}+r_n}(1)}= x_{r_n}(-1)\cdot
x_{r_n+2r_{n-1}}(-1),
\end{equation*}
 i.~e., $x_{n-1}=x_{r_{n-1}+r_n}(1)\cdot x\cdot x_{r_{n-1}+r_n}(1)$. Clearly
$x_{n-1}$ commutes with all $v_i$ for $k\le i\le n-3$. We construct $x_{n-3},\ldots,
x_k$ in the same way. Thus we obtain that $v$ and $v^{-1}$ are conjugate in $G$.
Again by \cite[Theorem~1.2 and~1.4]{Seitz} we have that all non-Abelian composition
factors of $C_G(v)$ are simple groups of Lie type of order less, than $Cmin$.
Therefore, $C_G(v)$ satisfies ($*$), a contradiction with Lemma \ref{power}. Hence
every $O_2(K_G)$-orbit $\Omega$ of $\Sym(I)$ with $O_3(K_G)\cap\langle X_r\vert
r\in\Omega\rangle\not=\{e\}$ has length $1$. Therefore, for some $i$, $O_3(K_G)\cap
X_{2e_i}\not=\{e\}$ is normal in $K$. Without lost we may assume that $i=1$ and
$X_{2e_1}=X_{r_0}$, where $r_0$ is the highest root of $\Phi$. Up to conjugation in
$G$ we may assume that $x_{r_0}(1)\leq Z(K)$. Thus $K_G\leq C_G(x_{r_0}(1))=P$,
where $P$ is the parabolic subgroup obtaining by removing $r_1$ from the set of
fundamental roots. The Levi factor $L$ of $P$ is known to have the following structure
$C_{n-1}(q)\ast S$, where $S=Z(L)$. Now, up to conjugation in $P$, $O_2(K_G)$ is the
Sylow $2$-subgroup of $C_{n-1}(q)$ and $KSO_3(P)/SO_3(P)$ is a Carter subgroup of
$C_{n-1}(q)\leftthreetimes \langle\zeta\rangle$. Since $t$ is odd we have that
$O_2(K_G)$ is a Sylow $2$-subgroup of $C_{n-1}(q)$ and Lemma
\ref{Syl2InCentrOfFieldAut} implies that $O_2(K_G)$ is a Sylow $2$-subgroup of
$(C_{n-1}(q))_\zeta=C_{n-1}(3)$. But, by Lemma \ref{Norm2Syl},
$N_{C_{n-1}(3)}(O_2(K_G))\not= O_2(K_G) C_{C_{n-1}(3)}(O_2(K_G))$, hence, by Lemma
\ref{CritSyl2Carter}, $K$ is not a Carter subgroup of $C_{n-1}(q)\leftthreetimes
\langle\zeta\rangle$. This final contradiction completes the proof.
\end{proof}

\section{Groups with triality automorphism}

\begin{ttt}\label{CarterTriality}
Let  $G$ be either
$D_4(q)$, or $ {}^3D_4(q^3)$. Assume that $\tau$ is a graph automorphism
of $G$ of order $3$ (in case of $ {}^3D_4(q^3)$ this is an automorphism, which
has the set of stable points isomorphic to  $D_4(q)$). Denote by $A_1$ the
subgroup of $\Aut(G)$ generated by all automorphisms, except $\tau$. Let
$A\leq \Aut(G)$ be such that $A\not\leq A_1$,  and  let $K$ be
a Carter subgroup of $A$.  Assume also that $\vert G\vert\le Cmin$.

Then  $\tau\in K$, $K\cap A_1$ is a Sylow $2$-subgroup of
$C_A(\tau)\simeq \Gamma G_2(q)$ and $\Gamma G_2(q)/G_2(q)$ is a $2$-group. A
Sylow $2$-subgroup $S$ of $\Gamma G_2(q)$ satisfies $N_{\Gamma
G_2(q)}(S)=SC_{\Gamma G_2(q)}(S)$ (i.~e. satisfies Lemma
{\em \ref{CritSyl2Carter}}) if either $q$ is odd,  or $q$ is even, $q=2^{2^t}$,
and~${\vert \Gamma G_2(q):G_2(q)\vert=2^t}$.
\end{ttt}

\begin{proof} In view of \cite[Theorem~1.2(vi)]{TiepZal} we have that every
element of $G$ is conjugate to its inverse. By \cite{Ca5} and \cite[Theorems~1.2
and~1.4]{Seitz} we obtain that for any element $t\in G$ of odd prime order, all
non-Abelian composition factors of $C_G(t)$ are simple groups of Lie type of order
less, then $Cmin$. Thus, $C_A(t)$ satisfies ($*$) and Lemma \ref{power} implies that
$K_G=K\cap G$ is a $2$-group. Now Lemma \ref{ConjAutomorphisms} implies that all
field or graph-field automorphisms of odd order of $G$ are conjugate under $G$. Since the
centralizer of every field or graph-field automorphism in $G$ is a group of
Lie type of order less, than $Cmin$, we again may apply Lemma \ref{power} and obtain that
$A$ does not contain field automorphism of odd order. Therefore, $K$ contains an element
$s$ of order $3$ such that $\langle s\rangle\cap A_1=\{e\}$, $G\leftthreetimes\langle
s\rangle=G\leftthreetimes\langle \tau\rangle$, and $K\cap A_1$ is a
$2$-group.

There exists two non-conjugate cyclic subgroups of order $3$: $\langle \tau\rangle$
and $\langle x\rangle$ such that $\langle\tau\rangle\cap A_1=\langle x\rangle\cap
A_1=\{e\}$ and $G\leftthreetimes\langle s\rangle=G\leftthreetimes\langle \tau\rangle$ (see
\cite[(9-1)]{GorLyo}). Hence, either $\tau\in K$, or $x\in K$.
Assume  that $q\not=3^t$. In the first case we obtain the statement of the theorem,
in the second case we have that $K\leq \langle K,C_G(x)\rangle$. But $C_G(x)\simeq\P
GL_3^\varepsilon(q)$, where $q\equiv \varepsilon=\pm1\pmod 3$ and $\P
GL_3^{+1}(q)=\P GL_3(q)$, $\P GL_3^{-1}(q)=\P GU_3(q)$. In view of \cite{DT2}, we
have that $K$ contains a Sylow $2$-subgroup of $\P GL_3^\varepsilon(q)$. Since the
normalizer of a Sylow $2$-subgroup of $\P GL_3^\varepsilon(q)$ is nilpotent, we have
that $K$ contains the normalizer of a Sylow $2$-subgroup of $\P
GL_3^\varepsilon(q)$. But, under our conditions, this normalizer contains a cyclic
subgroup $q-\varepsilon$ and $3$ divides $q-\varepsilon$. Therefore, $K\cap G$
contains an element of odd order, a contradiction. Thus this second case is
impossible.

Assume now that $q=3^t$. Then $C_G(\tau)\simeq G_2(q)$ and we obtain the theorem. In
the second case $C_G(x)\simeq SL_2(q)\rightthreetimes U$, where $U$ is a $3$-group
and $Z(C_G(x))\cap U\not=\{e\}$, a contradiction with Lemma~\ref{power}.
\end{proof}

\section{Carter subgroups in semilinear groups of Lie type}

\begin{ttt}\label{CarterSemilinear}
Let $G$ be a  finite group of Lie type  ($G$ is not necessary simple) over a field
of characteristic $p$  and $\overline{G}$, $\sigma$ are chosen so that
$O^{p'}(\overline{G}_\sigma)\leq G\leq\overline{G}_\sigma$. Assume also that $G\not\simeq
{}^3D_4(q^3)$. Choose a subgroup $A$
of $\Aut (O^{p'}(\ov{G}_\sigma))$ containing $G$ and assume that $A$ is not
contained in subgroup $A_1$ defined in Theorem \ref{CarterTriality}, if $G=D_4(q)$. Let
$K$ be a Carter subgroup of $A$
and assume that $A=KG$. Assume also that~$\vert G\vert\leq Cmin$.

Then exactly one of the following statements holds:

\begin{itemize}
 \item[{\em (1)}] $\Gamma G={}^2A_2(2^{2t})\leftthreetimes
\langle\zeta\rangle$, $\vert \zeta\vert=t$ is odd, $K\cap {}^2A_2(2^{2t})$ has
order~$2\cdot 3$ and is isomorphic to a Carter subgroup of~${^2A_2(2^2)}$.
\item[{\em (2)}] $G$ is defined over $GF(2^t)$, $\zeta$ is a field automorphism, $\vert\zeta\vert=t$, and
$K=S\leftthreetimes\langle\zeta\rangle$, where $S$ is a Sylow $2$-subgroup
 of~$G_{\zeta_{2'}}$.
\item[{\em (3)}] $G\simeq \P SL_2(3^{t})$, $\vert\zeta\vert=t$ is odd, and
$K=S\leftthreetimes\langle\zeta\rangle$, where $S$ is a   Sylow
$3$-subgroup of~$G_{\zeta_{3'}}$.
\item[{\em (4)}] $\Gamma G={}^2G_2(3^{2n+1})\leftthreetimes\langle\zeta\rangle$,
$\vert\zeta\vert=2n+1$, and $K\cap {}^2G_2(3^{2n+1})=S\times P$, where
$S$ is of order $2$ and~${\vert P\vert =3^{\vert\zeta\vert_3}}$.
\item[{\em (5)}] $p$ does not divide $\vert K\cap G\vert$ and $K$  contains  a
Sylow $2$-subgroup of~$A$. Note that in view of Lemmas \ref{InhBy2-ext},
\ref{Syl2InCentrOfFieldAut} and \ref{ConjAutomorphisms} a Sylow subgroup $S$ of $A$
satisfies $N_A(S)=S C_A(S)$ if $N_G(S\cap G)=(S\cap G) C_G(S\cap G)$.
\end{itemize}
In particular, Carter subgroups of $A$ are conjugate.
\end{ttt}

Note that after prooving this theorem we can do not demand that $A=KG$
and~${\vert A\vert\le Cmin}$.

We shall prove this theorem in the following way. If $G=C_n(q)$, then the
theorem follows from Theorem \ref{sympcarter}. If $G=A$ then the theorem follows
from \cite{PreTamVdo} and \cite{TamVdo}. Thus we may assume that $G\not=
C_n(q)$ and $G\not=A$. Then first of all, it is
possible to check by direct computation, that $G$ and $K$ can satisfy
precisely one of the statements of the theorem, i.~e., if (3) holds, and
$S$ is a Sylow $2$-subgroup of $A$, then, by Lemma \ref{Norm2Syl},
$N_A(S)\not=SC_A(S)$, hence ($5$) is not true. If ($4$) holds then by
\cite{war}, a Sylow $2$-subgroup $S$ of $A$  does not satisfy $N_A(S)=SC_A(S)$,
so by Lemma \ref{CritSyl2Carter}, statement (5) of the theorem is not true.

Assume that the theorem is false and $A$ is a counterexample with~$\vert G\vert$
minimal. Clearly this implies that $Z(G)$ is trivial. Let $K$ be a Carter subgroup
of $A$. First we show that if $p$ divides $\vert K\vert$, then one of (1)--(4) of
the theorem holds. Then we show that if $p$ does not divide $\vert K\vert$, then $K$
contains a Sylow 2-subgroup of $A$. Since both of these steps are quite complicated,
we divide them into two sections. Note also that by \cite{Ca5}, for every semisimple
element $t\in G$, all non-Abelian composition factors of $C_G(t)$, hence of $C_A(t)$
are simple groups of Lie type of order less, than $Cmin$. Hence $C_A(t)$ satisfies
($*$). In order to apply Lemma \ref{power} we shall use this fact without future
references. Further, for every non-Abelian composition factor $S\simeq D_4(q)$ of a
reductive
subgroup of maximal rank $R$ of $G$, we have that $\Aut_A(S)$ is contained in subgroup
$A_1\leq\Aut(S)$, i.~e. satisfies conditions of the theorem.

\section{Carter subgroups of order divisible by the characteristic}

Denote $K\cap G$ by $K_G$.
If $A$ contains a graph automorphism $\tau$ of $O^{p'}(\ov{G}_\sigma)$, then
every semisimple element of odd order is conjugate to its inverse in
$A$ (cf. Lemma \ref{ConjInverseInGraph}). By Lemma \ref{power} we obtain that
$\vert K_G\vert$ is divisible only by $2$ and $p$. If $p=2$, then we obtain that
$K_G$ is a $2$-group, it is contained in a proper $K$-invariant parabolic
subgroup $P$ of $G$  and by Lemma \ref{HomImageOfCarter} $KO_2(P)/O_2(P)$ is a
Carter subgroup of $KP/O_2(P)$. Since $K_G\leq O_2(P)$, it follows that
$\left(KO_2(P)/O_2(P)\right)\bigcap \left(P/O_2(P)\right)=\{1\}$. Hence $P$ is
 a Borel subgroup of $G$, otherwise we would have
$C_{P/O_2(P)}(KO_2(P)/O_2(P))\not=\{1\}$, a contradiction with the fact that $
KO_2(P)/O_2(P)$ is a Carter subgroup of $KP/O_2(P)$. Thus $P$ is a Borel
subgroup and the theorem follows from Lemma \ref{CartBorel}. Now if $p\not=
2$, then again $K_G$ is contained in a proper parabolic subgroup $P$ of $G$
such that $O_p(K_G)\leq O_p(P)$ and $O_2(K_G)\leq L$. Then Lemmas
\ref{Syl2InCentrOfFieldAut} and \ref{ConjAutomorphisms} implies that
$\Omega(H)\leq K$. Now Lemma \ref{centUH} implies that $O_p(K_G)\leq
C_U(\Omega(H))=\{e\}$. Therefore $K\cap G$ is a $2$-group. By Lemmas
\ref{Syl2InCentrOfFieldAut} and \ref{ConjAutomorphisms} every element  $x\in
A\setminus G$ of odd order such  that $\langle x\rangle\cap G=\{e\}$ centralizes
some Sylow $2$-subgroup of $G$. Hence $K$ contains a Sylow $2$-subgroup of $A$,
i.~e.,  $K$ satisfies (5) of the theorem.  Hence we may assume that~${A=\Gamma
G}$.

Recall that we are in the conditions of Theorem
\ref{CarterSemilinear}, $\Gamma G$ is supposed to be a counterexample
to the theorem with $\vert G\vert$ minimal and $K$ is a Carter
subgroup of $\Gamma G$ such that $p$ divides~$\vert K_G\vert$. We have
that $K=\langle \zeta^kg, K_G\rangle$. Since $\vert G\vert
\le Cmin$, Lemma \ref{HomImageOfCarter} implies that
 $KG/G$ is a Carter subgroup of $\Gamma G/G$. Therefore
$\vert\zeta^k\vert=\vert \zeta\vert$ and we may assume that~$k=1$.

In view of Lemma \ref{parabolic} there exists a proper $\sigma$- and
$\bar{x}$-invariant parabolic subgroup $\overline{P}$ of
$\overline{G}$ such that $O_p(K_G)\leq R_u(\overline{P})$ and $K_G\leq
\overline{P}$. In particular, $\overline{P}$ and
$\overline{P}^{\bar{\zeta}}$ are conjugate in $\overline{G}$. Let $\Phi$ be
the root system of $\ov{G}$ and
$\Pi$ be a set of fundamental roots of $\Phi$. In view of
\cite[Proposition~8.3.1]{Car1} $\overline{P}$ is conjugate to some
$\overline {P}_J=\overline{B}\cdot\overline{N}_J\cdot \ov{B}$, where $J$ is a
subset of $\Pi$ and $\overline{N}_J$ is a complete preimage of $W_J$
in $\overline{N}$ under the natural homomorphism
$\overline{N}/\overline{T}\rightarrow W$. Now $\overline{P}_J$ is
$\varphi$-invariant, hence
$\overline{P}_J^{\bar\zeta}=\overline{P}_J^{\bar\gamma^\varepsilon}$. Consider
the symmetry $\rho$  of the
Dynkin diagram of $\Phi$ corresponding to $\bar\gamma$. Let $\overline{J}$ be
the image of $J$ under
$\rho$. Clearly
$\overline{P}_J^{\bar\gamma}=\overline{P}_{\overline{J}}$. Since
$\overline{P}$ and $\overline{P}^{\bar\zeta}$ are conjugate in
$\overline{G}$ we obtain that $\overline{P}_J$ and
$\overline{P}_J^{\bar\zeta}$ are conjugate in $\overline{G}$. In view
of \cite[Theorem~8.3.3]{Car1} it follows that either $\varepsilon=0$,
or $J=\overline{J}$; i.~e.,  $\overline{P}_J$ is
$\bar\zeta$-invariant.

Now we have that $\overline{P}^{\bar{y}}=\overline{P}_J$. So
$\langle\bar{\zeta}g,\overline{P}\rangle^{\bar{y}}=\langle
(\bar{\zeta}g)^{\bar{y}},
\overline{P}_J\rangle$ and
$\overline{P}_J^{(\bar{\zeta}g)^{\bar{y}}}=\overline{P}_J$. It follows
$$(\bar{\zeta}g)^{\bar{y}}=\bar{y}^{-1}\bar{\zeta}g\bar{y}=
\bar{\zeta}\left(\bar{\zeta}^{-1}\bar{y}^{-1}\bar{\zeta}g\bar{y}\right)=
\bar{\zeta}\cdot h,$$ where
$h=\left(\bar{\zeta}^{-1}\bar{y}^{-1}\bar{\zeta}g\bar{y}\right)\in\overline{G}$.
Since $\overline{P}_J^{\bar\zeta}=\overline{P}_J=\overline{P}_J^{h^{-1}}$ we obtain
that $h\in N_{\overline{G}}(\overline{P}_J)$. By \cite[Theorem~8.3.3]{Car1},
$N_{\overline{G}}(\overline{P}_J)=\overline{P}_J$, thus
$\langle\bar{\zeta}g,\overline{P}\rangle^y=\langle \bar\zeta,\overline{P}_J\rangle$.
Now both $\overline{P}$ and $\overline{P}_J$ are $\sigma$-invariant. Hence
$\bar{y}\sigma(\bar{y}^{-1})\in N_{\overline{G}}(\overline{P})=\overline{P}$.
Therefore, by Lang-Steinberg Theorem \cite[Theorem~10.1]{Ste}  we may assume that
$\bar{y}=\sigma(\bar{y})$, i.~e. $\bar{y}\in\overline{G}_\sigma$. Since
$\overline{G}_\sigma=\overline{T}_\sigma\cdot O^{p'}(\overline{G}_\sigma)$ and
$\overline{T}\leq \overline{P}_J$, then we may assume that $\bar{y}\in
O^{p'}(\overline{G}_\sigma)$. Thus, up to conjugation in $G$, we may assume that
$\overline{K}\leq \langle \bar\zeta,\overline{P}_J\rangle=
\overline{P}_J\leftthreetimes\langle\bar\zeta\rangle$ and $K\leq (\overline{P}_J\cap
G)\leftthreetimes\langle \zeta\rangle=P_J\leftthreetimes\langle\zeta\rangle$.
Further if $\overline{L}_J=\langle\overline{T}, \overline{X}_r\vert r\in J\cup
-J\rangle$, then $\overline{L}_J$ is a $\sigma$- and $\bar\zeta$-invariant Levi
factor of $\overline{P}_J$ and $L_J=\overline{L}_J\cap G$ is a $\zeta$-invariant
Levi factor of~$P_J$.  Lemma \ref{HomImageOfCarter} implies that
$KO_p(P_J)/O_p(P_J)=X$ is a Carter subgroup of $L_J\leftthreetimes\langle
\zeta\rangle$ and $K_1Z(L_J)/Z(L_J)=\widetilde{X}$ is a Carter subgroup
of~$(L_J\leftthreetimes \langle \zeta\rangle)/Z(L_J)$. Recall that $K=\langle \zeta
g,K_G\rangle$, hence, if $v$ and $\tilde{v}$ are the images of $g$ under the natural
homomorphisms $\omega:P_J\leftthreetimes\langle\zeta\rangle\rightarrow
L_J\leftthreetimes\langle\zeta\rangle\simeq(P_J
\leftthreetimes\langle\zeta\rangle)/O_p(P_J)$ and
$\tilde\omega:P_J\leftthreetimes\langle \zeta\rangle\rightarrow
(P_J\leftthreetimes\langle\zeta\rangle)/Z(L_J)O_p(P_J)\simeq
\left(L_J\leftthreetimes\langle \zeta\rangle\right)/Z(L_J)$, then $X=\langle\zeta
v,K_G^{\omega}\rangle$ and~$\widetilde{X}=\langle\zeta
\tilde{v},K_G^{\tilde\omega}\rangle$. Note that $O_p(P)$ and $Z(L_J)$ are
characteristic subgroups of $P$ and $L_J$ respectively, hence we may consider
$\zeta$ as an automorphism of $L_J\simeq P/O_p(P)$ and $\widetilde{L}=L_J/Z(L_J)$.
Note also that all non-Abelian composition factors of $P$ are simple groups of Lie
type of order less than $Cmin$, hence $P\leftthreetimes\langle\zeta\rangle$
satisfies ($*$). Thus we may apply Lemma \ref{HomImageOfCarter}
to~${\widetilde{L}\leftthreetimes\langle\zeta\rangle}$.

If $P_J$ is a Borel subgroup of $G$, we use Lemma \ref{CartBorel}. So we may assume
that $L_J\not=Z(L_J)$, i.~e., that $P_J$ is not a Borel subgroup of $G$. Then
$L_J=H(G_1\ast\ldots\ast G_k)$, where $G_i$ are subsystem subgroups of $G$ and $H$
is a Cartan subgroup of $G$. Let $\zeta=\zeta_2\cdot\zeta_{2'}$ be the product of
$2$- and $2'$-parts of $\zeta$. Now $\zeta_{2'}=\varphi^k$, for some $k$, is a field
automorphism (recall that we do not consider the triality automorphism) and it
normalizes every $G_i$, since $\varphi$ normalizes every $G_i$. Moreover, in view of
Lemma \ref{Syl2InCentrOfFieldAut}, we have that $\zeta_{2'}$ centralizes a Sylow
$2$-subgroup of $H$. In particular, it centralizes a Sylow $2$-subgroup of
$Z(L_J)\leq H$. Therefore, any element of odd order of $\langle
\zeta\rangle\rightthreetimes L_J$ centralizes a Sylow $2$-subgroup of~$Z(L_J)$.

Now $\widetilde{L}=(\P G_1\times\ldots\times \P
G_k)\widetilde{H}$, where
$\widetilde{H}=H^{\omega_2}$ and $\P G_1,\ldots,\P G_k$ are canonical finite
groups of
Lie type with trivial centre. Denote by $M_i=C_{\widetilde{L}}(G_i)$, clearly
$M_i=(\P G_1\times
\ldots\times \P G_{i-1}\times \P
G_{i+1}\times\ldots\times \P G_k)C_{\widetilde{H}}(G_i)$; denote by $L_i$ the
factor
group
$\widetilde{L}/M_i$ and by $\pi_i$ corresponding natural
homomorphism. Then $L_i$ is a finite group of Lie type and $\P G_i\leq
L_i\leq \widehat{\P G}_i$. Denote by $M_{i,j}=C_{\widetilde{L}}(\P G_i\times \P
G_j)$, then $M_{i,j}=(\P G_1\times
\ldots\times \P G_{i-1}\times \P
G_{i+1}\times\ldots\times \P G_{j-1}\times \P G_{j+1}\times\ldots\times\P
G_k)C_{\widetilde{H}}(\P G_i\times\P G_j)$; denote by $\pi_{i,j}$
corresponding natural homomorphism $\widetilde{L}\rightarrow
\widetilde{L}/M_{i,j}$. If $M_i$ (respectively $M_{i,j}$) is $\zeta$-invariant,
then $M_i$ (resp. $M_{i,j}$) is normal in
$\widetilde{L}\leftthreetimes\langle\zeta\rangle$ and we denote by $\delta_i$
(resp. $\delta_{i,j}$) the natural homomorphism
$\delta_i:\widetilde{L}\leftthreetimes\langle\zeta\rangle\rightarrow
\left(\widetilde{L}
\leftthreetimes\langle\zeta\rangle\right)/M_i$
($\delta_{i,j}:\widetilde{L}\leftthreetimes\langle\zeta\rangle\rightarrow
\left(\widetilde{L }
\leftthreetimes\langle\zeta\rangle\right)/M_{i,j}$).

Now consider $\zeta_2$. Since $\zeta_2^2$ is a field automorphism, there can be
two cases: either $\zeta_2$
normalizes $\P G_i$, or $\zeta^2_2$ normalizes $\P G_i$ and
$\P G_i^{\zeta_2}=\P G_j$ for some $j\not=i$. Consider these two cases
separately.

Case 1, $\zeta_2$ normalizes $\P G_i$. Then $\zeta_2$ normalizes
$M_i$, and Lemma \ref{HomImageOfCarter} implies that
$\widetilde{X}^{\delta_i}=K_i$ is a Carter subgroup of
$L_i\leftthreetimes\langle
\zeta\rangle$. Since $L_i\leftthreetimes\langle
\zeta\rangle$ is
a semilinear group of Lie type satisfying the conditions of Theorem
\ref{CarterSemilinear},  $\vert L_i \vert<\vert G\vert$,
and $p$ does not divide $\vert K_i\vert$, we have that  $K_i$ contains a
Sylow $2$-subgroup $S_i$ of $L_i\leftthreetimes\langle
\zeta\rangle$ (in particular, $p=2$)
and, by Lemma \ref{CritSyl2Carter}, $N_{L_i\leftthreetimes\langle
\zeta\rangle}(S_i)=S_i C_{L_i\leftthreetimes\langle
\zeta\rangle}(S_i)=S_i\times
O(N_{L_i\leftthreetimes\langle
\zeta\rangle}(S_i))$.

Case 2, $\zeta_2^2$ normalizes $\P G_i$ and $\P G_i^{\zeta_2}=\P G_j$. Then
$M_{i,j}$ is normal in $\widetilde{L}\leftthreetimes\langle\zeta\rangle$. We want to
show that if $S_{i,j}$ is a Sylow $2$-subgroup of
$(\widetilde{L})^{\pi_{i,j}}\leftthreetimes\langle\zeta\rangle$, then
$N_{(\widetilde{L})^{\pi_{i,j}}\leftthreetimes\langle\zeta\rangle}(S_{i,j})=S_{i ,j}
C_{(\widetilde{L})^{\pi_{i,j}}\leftthreetimes\langle\zeta\rangle}
(S_{i,j})=S_{i,j}\times
O(N_{(\widetilde{L})^{\pi_{i,j}}\leftthreetimes\langle\zeta\rangle}(S_i))$. Since
$M_{i,j}$ is a normal subgroup of $\widetilde{L}\leftthreetimes\langle\zeta\rangle$,
then,  by Lemma \ref{HomImageOfCarter}, $(\widetilde{X})^{\delta_{i,j}}$ is a Carter
subgroup of $(\widetilde{L}\leftthreetimes\langle\zeta\rangle)^{\delta_{i,j}}$. Thus
we may assume that $\widetilde{L}=\P G_1\times \P G_2$ and $(\P G_1)^\zeta=\P G_2$.
Now we are in the condition of Lemma \ref{CarterInDirectProduct}, namely, we have a
finite group $\widetilde{G}=\widetilde{X}(\P G_1\times \P G_2)$, where $\P
G_1\simeq\P G_2$ has trivial centre. Then  $\Aut_{\widetilde{X}}(\P G_1)$ is a
Carter subgroup of $\Aut_{\widetilde{G}}(\P G_1)$. Now $\P G_1$ is a canonical finite
group of
Lie type and $\P G_1\leq \Aut_{\widetilde{G}}(\P G_1)\leq\Aut (\P G_1)$, i.~e.,
$\Aut_{\widetilde{G}}(\P G_1)$  satisfies conditions of Theorem
\ref{CarterSemilinear} and $\widetilde{X}\cap (\P G_1\times\P G_2)$ is not divisible
by the characteristic. By induction, $\Aut_{\widetilde{X}}(\P G_1)$ contains a Sylow
$2$-subgroup of $\Aut_{\widetilde{G}}(\P G_1)$ (in particular, $p=2$). The same arguments
show that
$\Aut_{\widetilde{X}}(\P G_2)$ contains a Sylow $2$-subgroup of
$\Aut_{\widetilde{G}}(\P G_2)$. Now consider $N_{\widetilde{G}}(\P
G_1)=N_{\widetilde{G}}(\P G_2)$. Let $S$ be a Sylow $2$-subgroup of
$N_{\widetilde{G}}(\P G_1)$. Since $C_{\widetilde{G}}(\P G_1)\cap
C_{\widetilde{G}}(\P G_2)=\{1\}$, Lemma \ref{Syl2centrcomposit} implies that
$N_{N_{\widetilde{G}}(\P G_1)}(S)=S C_{N_{\widetilde{G}}(\P G_1)}(S)$. Now $\vert
\widetilde{G}:N_{\widetilde{G}}(\P G_1)\vert=2$, thus, by Lemma \ref{InhBy2-ext}, we
obtain that $N_{\widetilde{G}}(\widetilde{S})=\widetilde{S}C_{N_{\widetilde{G}}(\P
G_1)}(\widetilde{S})$ for a Sylow $2$-subgroup $\widetilde{S}$ of $\widetilde{G}$.
Hence, if $S_{i,j}$ is a Sylow $2$-subgroup of $(\P G_i \times \P
G_j)\leftthreetimes\langle\zeta\rangle$, then $$N_{(\P G_i\times \P
G_j)\leftthreetimes\langle \zeta\rangle}(S_{i,j})=S_{i,j} C_{(\P G_i\times \P
G_j)\leftthreetimes\langle \zeta\rangle}(S_{i,j})=S_{i,j}\times O(N_{(\P G_i\times
\P G_j)\leftthreetimes\langle \zeta\rangle}(S_{i,j})).$$

Now we show that for a Sylow $2$-subgroup $S$ of $L_J\leftthreetimes
\langle\zeta\rangle$ we have that $N_{L_J\leftthreetimes
\langle\zeta\rangle}(S)=SC_{L_J\leftthreetimes \langle\zeta\rangle}(S)$. Since
$\widetilde{L}\not=\{e\}$, then, as we noted above, $p\not=2$.
Consider an element $x\in N_{L_J\leftthreetimes \langle\zeta\rangle}(S)$ of odd
order. We  need to prove that $x$ centralizes $S$. As we noted above, any
element of odd order of
$L_J\leftthreetimes\langle\zeta\rangle$ centralizes $S\cap Z(L_J)$, hence if $x$
centralizes $S/(S\cap Z(L_J))$, then $x$ centralizes $S$.  Now either $M_i$ is
normal in
$\widetilde{L}\leftthreetimes
\langle\zeta\rangle$, or $M_{i,j}$ is normal in
$\widetilde{L}\leftthreetimes\langle\zeta\rangle$ and
$\left(\cap_iM_i\right)\bigcap\left(\cap_{i,j}M_{i,j}\right)=\{1\}$. Moreover,
as we proved above $x^{\delta_i}$
centralizes $SM_i/M_i$, and $x^{\delta_{i,j}}$ centralizes $S M_{i,j}/M_{i,j}$.
By Lemma \ref{Syl2centrcomposit} we obtain that $x$
centralizes~$S$.

Thus $N_{L_J\leftthreetimes \langle\zeta\rangle}(S)=SC_{L_J\leftthreetimes
\langle\zeta\rangle}(S)$, by Lemma \ref{CritSyl2Carter} there exists a Carter
subgroup $F$ of $L_J\leftthreetimes \langle\zeta\rangle$ containing $S$. Since
$L_J\leftthreetimes \langle \zeta\rangle$ satisfies ($*$), Theorem
\ref{maininduct} implies that $K$ and $F$ are conjugate, i.~e. $K$ contains a
Sylow $2$-subgroup of $L_J\leftthreetimes \langle\zeta\rangle$. In particular,
$\Omega(H)\leq K$ and $\Omega(H)$ centralizes $K_p=K\cap O_p(P_J)\not=\{e\}$. A
contradiction with Lemma~\ref{centUH}.

\section{Carter subgroups of order not divisible by the characteristic}

Again we are in the conditions of Theorem \ref{CarterSemilinear}. If
$A$ contains a graph automorphism $\gamma$, or $\Phi(\ov{G})\not=A_n,D_{2n+1
}, E_6$, then by Lemma
\ref{ConjInverseInGraph} and \cite[Lemma~2.2]{TamVdo} we obtain that every
semisimple element of
$G$ is conjugate to its inverse. By Lemma \ref{power} we have that
$K_G=K\cap G$ is a $2$-group. By Lemmas \ref{Syl2InCentrOfFieldAut} and
\ref{ConjAutomorphisms} we obtain
that a $2'$-part of $K$ centralizes a Sylow $2$-subgroup of $G$, hence
$K_G$ is a Sylow $2$-subgroup of $G$ and $K$ contains a Sylow $2$-subgroup of
$A$. Thus Theorem
\ref{CarterSemilinear} is true in this case. So we may assume that
$A=\Gamma G=\langle\zeta,G\rangle$  is a semilinear group of Lie type,
$K=\langle \zeta^kg,K_G\rangle$ is a Carter
subgroup of $\Gamma G$, and $\Phi(\ov{G})\in\{A_n, D_{2n+1}, E_6\}$. Like in the
previous section we may assume
that~$k=1$. Since $G_\zeta$ is nontrivial we have that $K_G$ is
nontrivial also. Therefore $Z(K)\cap K_G$ is nontrivial. Consider an element
$x\in Z(K)\cap K_G$ of prime order. Then $K\in C_{\Gamma G}(x)=\langle \zeta
g,C_G(x)\rangle$. Now $C_{\overline{G}}(x)^0=\overline{C}$ is a
connected $\sigma$-stable reductive subgroup of maximal rank of
$\overline{G}$. Moreover $\overline{C}$ is a characteristic subgroup
of $C_{\overline{G}}(x)$ and $C_{\ov{G}}(x)/\ov{C}$ is isomorphic to a
subgroup of $\Delta$ (see \cite[F, \S4 and
Proposition~5]{BorCar}). Thus $K$ is contained in $\langle K,
C\rangle$, where $C=\ov{C}\cap G$. Moreover $C=\ov{C}\cap
G=T(G_1\ast\ldots\ast G_m\ast S)$ is normal in $C_{\Gamma
G}(x)$ and $K_GC/C$ is isomorphic to a subgroup of $\Delta$.  There
can be two cases:

$m=0$ and $\vert x\vert>2$, i.~e., $C=T=S$ is a maximal torus. Then $\ov{T}$ is
$\bar\zeta g$-stable. In view of Lemma \ref{NormOfRegularElementIsNotCentr} we
obtain that there exists $h\in O^{p'}(\ov{G}_\sigma)$ normalizing, but not
centralizing $\langle x\rangle$. Therefore $x$ is conjugate to its nontrivial power,
a contradiction with Lemma~\ref{power}.

Either $m\ge1$ or $\vert x\vert=2$. Assume first that $m\ge1$ and $\vert x\vert
> 2$. Then $Z(C)=S$ and $G_1\ast\ldots\ast G_m$ are normal subgroups of
$\langle K, C\rangle$. Hence we may consider $\widetilde{G}=\langle
K,G_1\ast\ldots\ast G_m\ast S\rangle/S$. Then $\widetilde{G}=\widetilde{K}(\P
G_1\times\ldots\times \P G_m)$, where $\widetilde{K}=KZ(C)/Z(C)$ is a Carter
subgroup of $\widetilde{G}$ (cf. Lemma \ref{HomImageOfCarter}) and $Z(\P G_i)$ is
trivial. Now $\widetilde{K}$ acts by conjugation on $\{\P G_1,\ldots, \P G_ m\}$ and
without lost we may assume that $\{\P G_1,\ldots,\P G_m\}$ is a
$\widetilde{K}$-orbit. Thus we are in the condition of Lemma
\ref{CarterInDirectProduct} and $\Aut_{\widetilde{K}}(\P G_1)$ is a Carter subgroup
of $\Aut_{\widetilde{G}}(\P G_1)$. Moreover $\vert\widetilde{K}\cap \P
G_1\times\ldots\times\P G_m\vert$ is not divisible by the characteristic. By
induction we have that $\Aut_{\widetilde{K}}(\P G_1)$ contains a Sylow $2$-subgroup
of $\Aut_{\widetilde{G}}(\P G_1)$, hence the characteristic is odd and $\vert K\cap
G\vert$ is divisible by $2$. Therefore we may assume $x\in Z(K)\cap K_G$ to be an
involution. Furthermore, Lemma \ref{CentrOfInvolution} implies that $\vert
K_G:(K_G\cap C)\vert$ divides $\vert C_G(x)/C\vert$ and $C_G(x)/C$ is a $2$-group.

Thus $x$ is an involution. By Lemma \ref{InvolutionsAndTori}  we have that
every involution of
$G$ is contained in a maximal torus $T$ such that $N(G,T)/T\simeq W$, where $W$
is the Weyl group of $\ov{G}$. In particular, we may assume that
$\ov{C}$ is generated by $\overline{T}$-invariant root
subgroups and that $\zeta_{2'}$ centralizes $C_G(x)/C$ (by
\cite[Proposition~2]{Ca5} $N(G,C)/C\simeq N_W(W_1)/W_1$, where $W=N(G,T)/T$ and
$W_1=N(C,T)/T$). In view of Lemma
\ref{Syl2InCentrOfFieldAut}, $\zeta_{2'}$ normalizes $T$, each of $G_i$ and
centralizes some Sylow $2$-subgroups of $T$ and of each of $G_i$. We can
write $\zeta g=\zeta_2 g_1\cdot \zeta_{2'} g_2$, where $\zeta_2g_1$ is
a $2$-part and $\zeta_{2'}g_2$ is a $2'$-part of $\zeta g$. Since $\zeta_{2'}$
centralizes $C_G(x)/C$ and since $C_G(x)/C$ is a $2$-group, we obtain that
$g_2\in C$. Hence
$\zeta_{2'}g_2$ normalizes $T$ and each of $G_i$. As we noted above,
$\zeta_{2'}$ centralizes a Sylow $2$-subgroup of $T$, hence, it centralizes a
Sylow $2$-subgroup of $Z(C)\leq T$. Thus we have that every element of odd order
of $\langle K,C_G(x)\rangle$
centralizes the Sylow $2$-subgroup of $Z(C)$.

Now consider $\widetilde{G}=\langle K,C\rangle/Z(C)$. Then
$\widetilde{G}=\widetilde{K}(\P G_1\times\ldots\times\P G_m)$ (possibly $m=0$),
where $\widetilde{K}=KZ(C)/Z(C)$ is a Carter subgroup of $\widetilde{G}$ (cf. Lemma
\ref{HomImageOfCarter}) and for all $i$, ${Z(\P G_i)=1}$. By Lemma
\ref{CarterInDirectProduct} we have that $\Aut_{\widetilde{K}}(\P G_1)$ is a Carter
subgroup of $\Aut_{\widetilde{G}}(\P G_1)$. Since $\P G_1$ is a finite group of Lie
type satisfying Theorem \ref{CarterSemilinear}, by induction we obtain that
$\Aut_{\widetilde{K}}(\P G_1)$ contains a Sylow $2$-subgroup of
$\Aut_{\widetilde{G}}(\P G_1)$. Similarly we have that $\Aut_{\widetilde{K}}(\P
G_i)$ contains a Sylow $2$-subgroup of $\Aut_{\widetilde{G}}(\P G_i)$ for all $i$.
Let $S$ be a Sylow $2$-subgroup of $N_{\widetilde{G}}(\P G_1)$. Since
$C_{\widetilde{G}}(\P G_1\times\ldots\times \P G_m)=\{1\}$, Lemma
\ref{Syl2centrcomposit} implies that $N_{N_{\widetilde{G}}(\P G_1)}(S)=S
C_{N_{\widetilde{G}}(\P G_1)}(S)$. Now $\vert \widetilde{G}:N_{\widetilde{G}}(\P
G_1)\vert=2^t$, thus, by Lemma \ref{InhBy2-ext}, we obtain that
$N_{\widetilde{G}}(\widetilde{S})=\widetilde{S}C_{N_{\widetilde{G}}(\P
G_1)}(\widetilde{S})$ for a Sylow $2$-subgroup $\widetilde{S}$ of $\widetilde{G}$.
Since $\vert\P G_i\vert< Cmin$, then $\widetilde{G}$ and $\langle K,C\rangle$
satisfy ($*$). By Lemma \ref{CritSyl2Carter} we obtain that there exists a Carter
subgroup  $M$ of $\langle K,C\rangle$ that contains a Sylow $2$-subgroup of $\langle
K,C\rangle$. By Theorem \ref{maininduct}, subgroups $M$ and $K$ are conjugate in
$\langle K,C\rangle$, thus  $K$ contains a Sylow $2$-subgroup $ Q$ of $\langle
K,C_G(x)\rangle$.

Let $S_1$ be a Sylow $2$-subgroup of $\Gamma G$ containing $Q$ and let $t\in
Z(S_1)\cap G$. Then $t\in C_G(x)$, hence, $t\in Z(Q)$ and $t\in Z(K)$. Thus we
may substitute
$x$ by $t$ in arguments above and obtain that $Q=S_1$, i.~e., $K$ contains a
Sylow $2$-subgroup of~$\Gamma G$.

\section{Carter subgroups of finite groups are conjugate}

In order to state the following theorem without using the classification of
finite simple groups, we give the following definition. A finite group is said
to be a {\em $K$-group} if all its non-Abelian composition factors are known
simple groups.

\begin{ttt}\label{main}
Let $G$ be a finite $K$-group. Then Carter subgroups of $G$ are conjugate.
\end{ttt}

\begin{proof}
By \cite[Theorems~3.3--3.5]{TamVdo}, \cite[Theorem~1.1]{PreTamVdo},
\cite[Table]{Vdo}; and Theorems \ref{sympcarter}, \ref{CarterTriality}, and
\ref{CarterSemilinear} from the present paper we obtain that $G$ satisfies ($*$).
Hence, by Theorem \ref{maininduct}, Carter subgroups of $G$ are conjugate.
\end{proof}

\end{document}